\documentclass[12pt]{article}
\usepackage{e-jc}

\usepackage[T1]{fontenc}

\usepackage{amsmath,amssymb}
\usepackage{mathtools}
\usepackage{enumitem}
\usepackage{graphicx}
\usepackage[colorlinks=true,citecolor=black,linkcolor=black,urlcolor=blue]{hyperref}

\newcommand{\disc}{\mathrm{Disc}}
\newcommand{\qdisc}{\mathrm{QDisc}}

\newcommand{\tr}{\mathrm{tr}}
\newcommand{\diag}{\mathrm{diag}}

\newcommand{\ra}{\mathrm{rank}}
\newcommand{\cc}{\raisebox{1.5pt}{$,$}}
\newcommand{\vol}{\mathrm{vol}}
\newcommand{\rchi}{\raisebox{1.9pt}{$\chi$}}

\newcommand{\lrchi}{\raisebox{.87pt}{\scalebox{.83}{$\chi$}}}

\newcommand{\norm}[2]{\left\|{#1}\right\|_{{#2}}}
\newcommand{\cl}[1]{\mathcal{#1}}
\newcommand{\bb}[1]{\mathbb{#1}}
\newcommand{\wt}[1]{\widetilde{#1}}

\title{Quantum Discrepancy: A Non-Commutative Version of Combinatorial Discrepancy}

\author{Kasra Alishahi\\
\small Department of Mathematical Sciences\\[-0.8ex]
\small Sharif  University of Technology\\[-0.8ex] 
\small Tehran, Iran\\
\small\tt alishahi@sharif.ir\\
\and
Mohaddeseh Rajaee\qquad Ali Rajaei\\
\small Faculty of Mathematical Sciences\\[-0.8ex]
\small Tarbiat Modares University\\[-0.8ex]
\small Tehran, Iran\\
\small\tt \{mohaddeseh.rajaee,alirajaei\}@modares.ac.ir}
\date{}
\begin{document}

\maketitle


\begin{abstract}
  In this paper, we introduce a notion of quantum discrepancy, a non-commutative version of combinatorial discrepancy which is defined for projection systems, i.e. finite sets of orthogonal projections, as non-commutative counterparts of set systems. We show that besides its natural algebraic formulation, quantum discrepancy, when restricted to set systems, has a probabilistic interpretation in terms of determinantal processes. Determinantal processes are a family of point processes with a rich algebraic structure.  A common feature of this family is the local repulsive behavior of points. Alishahi and Zamani (2015) exploit this repelling property to construct low-discrepancy point configurations on the sphere. 
  
  We give an upper bound for quantum discrepancy in terms of $N$, the dimension of the space, and $M$, the size of the projection system, which is tight in a wide range of parameters $N$ and $M$. Then we investigate the relation of these two kinds of discrepancies, i.e. combinatorial and quantum, when restricted to set systems, and bound them in terms of each other.
\end{abstract}

\section{Introduction}
\subsection{Definition of Quantum Discrepancy}
\label{Sub: Definition of QuantDisc}
In 1964, Roth proved in~\cite{Roth} that for any blue-red coloring of $[N]\coloneqq\{1,2,\ldots,N\}$, there always exists an arithmetic progression, in which the difference between the number of red and blue points is $\Omega(N^{1/4})$. Roth's theorem was one of the first results in a field that was later named combinatorial discrepancy theory, which is also  related to, or has applications in diverse fields of mathematics and theoretical computer science, such as  Ramsey theory, hypergraph coloring, arithmetic structures, probabilistic and approximation algorithms, complexity theory, and data structure. For a general introduction to discrepancy theory, its relations to computer science, and recent results in this field see~\cite{Chazelle,Panorama,Matousek}.

If we forget the special structure of the family of arithmetic progressions, we can formulate the basic problem of this field in a more general form.
Assume that $\Omega$ is a finite set. A \emph{$2$-coloring} or more simply a \emph{coloring} of $\Omega$ can be modeled by a function $\rchi: \Omega \to \{\pm 1\}$. Assuming $\rchi(s)=1$ and $\rchi(s)=-1$ as $s$ is colored red and blue, respectively, the absolute value of $\rchi(S)\coloneqq\sum_{s\in S}\rchi(s)$ is a measure of the imbalance between the number of red and blue elements  in $S$ due to $\rchi$. Now, given a \emph{set system} $\cl S$, i.e. a subset of $2^{\Omega}$, the \emph{discrepancy} of  $\cl{S}$ is defined by
\begin{align}
	\label{eq:CombDisc}
	\disc(\cl{S})&=\min_{\lrchi:\,\text{coloring}}\max_{S\in\cl{S}}|\rchi(S)|
\end{align}
Note that since the structure of a set system remains unaffected if we change the labels of points in the ground set, from now on, we suppose that $\Omega=[N]$ for some $N\in\bb{N}$.

Since the combinatorial minimization in \eqref{eq:CombDisc} has exponentially many feasible values, exhaustive search is not a computationally efficient method for finding a \emph{low-discrepancy coloring}, i.e. a coloring such as $\rchi$ that makes $\max_{S\in\cl S}|\rchi(S)|$ small. Moreover, many of general upper bounds such as Spencer's bound, mentioned in the next subsection, are based on techniques that are non-constructive, in the sense that they prove the existence of low-discrepancy colorings without any guide to accessing them. Designing efficient algorithms for the construction of colorings with optimal or nearly optimal discrepancy is an active area of research in the theory of combinatorial discrepancy. In recent years, remarkable achievements have been attained in this area (e.g. see~\cite[Ch.\ 6]{Panorama} and~\cite{Rothvoss}).

Quantum discrepancy is a non-commutative version of combinatorial discrepancy, which is defined as follows. In what follows, $\cl{M}_N(\mathbb C)$ stands for the set of $N\times N$ matrices with complex entries.
\begin{definition}
	\label{Def: QuantDisc}
	Suppose $N$ is a natural number. A finite set $\cl{P}$ of orthogonal projections of $\bb{C}^N$, i.e. operators such as $P$ with the property that $P^2 = P=P^*$, is called a \emph{projection system} in $\bb{C}^N$. A Hermitian matrix $\rchi\in\cl{M}_N(\bb C)$ with eigenvalues in $\{\pm 1\}$ is named a \emph{quantum coloring} (a justification for this denomination is provided at the end of this subsection). We define the \emph{quantum discrepancy} of a projection system $\cl{P}$ to be
			\begin{align}
			\label{eq: NonCommutDef}
			\qdisc(\cl{P})=\min_{\substack{{\lrchi: \text{\,quantum}}\\{\text{coloring}}}}\max_{P\in\cl{P}}\left[\tr^2(\rchi P)+\tr\left(\rchi[\rchi,P]P\right)\right]^{\frac{1}{2}},
			\end{align}
			where $[A,B]\coloneqq AB-BA$ is the commutator of  $A$ and $B$.
\end{definition}

Quantum discrepancy is well-defined for each projection system since the expression $\tr^2(\rchi P)+\tr\left(\rchi[\rchi,P]P\right)$ is real and non-negative. Because  $\rchi$ and $P$ are Hermitian, $\tr(\rchi P)$ is real. Also, it will be seen in Subsection~\ref{SubSec: NormCommutForm} that $\tr\left(\rchi[\rchi,P]P\right)$ is real and non-negative (see Lemma~\ref{lem: NormTrace}).

To a set system $\cl{S}\subseteq2^{[N]}$ we can assign a projection system $\cl{P}_{\cl{S}}$ in a natural way: for $S\in\cl{S}$ we set $P_S$ to be the diagonal matrix that for each $i\in[N]$, its $i$-th diagonal entry is equal to $1$ if $i\in S$, and $0$ otherwise. Then, we put $\cl P_{\cl S}=\{P_S: S\in\cl S\}$. Moreover, any coloring $\rchi$ of $[N]$ can be considered as a diagonal matrix with diagonal entries $\pm 1$. Thus, projection systems and quantum colorings generalize set systems and colorings, respectively. In addition, for a set $S\subseteq[N]$, and a (diagonal) coloring $\rchi$, $\tr(\rchi P_S)=\rchi(S)$ and since $\rchi$ and $P_S$ are both diagonal, $[\rchi,P_S]=\mathbf{0}$. Hence, in this situation,
\begin{align*}
   \left[\tr^2(\rchi P_S)+\tr\left(\rchi[\rchi,P_S]P_S\right)\right]^{\frac{1}{2}}=|\tr(\rchi P_S)|=\left|\rchi(S)\right|\cdot
\end{align*} 
It is proved in Subsection~\ref{SubSec: NormCommutForm} that, for an orthogonal projection $P$ and a quantum coloring $\rchi$, it holds that
\begin{align*}
\tr\left(\rchi[\rchi,P]P\right)=\frac{1}{2}\norm{[\rchi,P]}{2}^2.
\end{align*} 
Therefore, the second term, i.e.  $\tr\left(\rchi[\rchi,P]P\right)$, has also a simple meaning. The reason for choosing the combination $\left[\tr^2(\rchi P)+\frac{1}{2}\norm{[\rchi,P]}{2}^2\right]^{1/2}$ which might look messy, lies in its probabilistic interpretation, 
when $\cl P$ is associated to a set system. In fact, for a quantum coloring $\rchi$ there is a determinantal point process $\mathfrak{X}$ on $[N]$ (see Definition~\ref{Def: DetPross}) 
such that for any $S\subseteq[N]$,
\begin{align}
\label{eq: ExpKer}
\left[\tr^2(\rchi P_S)+\tr\left(\rchi[\rchi,P_S]P_S\right)\right]^{\frac{1}{2}}
=\left[\bb E\left[ \left(2\mathfrak X(S)-|S|\right)^2\right]\right]^{\frac{1}{2}}. 
\end{align}
In this formula $\mathfrak X(S)$ represents $|\mathfrak{X}\cap S|$, i.e. the number of points of $\mathfrak X$ which belong to $S$. Therefore, if we assume that the set of red points of $[N]$ is determined randomly according to the distribution of $\mathfrak{X}$, the right-hand side of \eqref{eq: ExpKer} can be considered as the imbalance imposed by the process $\mathfrak X$ to $S$. This probabilistic approach to the  combinatorial discrepancy problem (including the poof of \eqref{eq: ExpKer}) is fully described in Subsection~\ref{Sub: ProbApp}.

We intend to emphasize the importance of the commutator operator in \eqref{eq: NonCommutDef}. The appearance of $[\rchi ,P]$ reveals the strong non-commutative essence of the quantum discrepancy. This non-commutativity is partially a result of extending the notion of the coloring. Thus, we have a non-commutative concept of discrepancy even for set systems. There is another aspect of non-commutativity which lies in $\cl{P}$. When all projections in $\cl{P}$ commute pairwise, they are simultaneously unitarily diagonalizable since they are Hermitian. Therefore, $\cl {P}$, possibly after a unitary change of the basis, corresponds to a set system, and because of the cyclic property of the trace function, $\qdisc(\cl{P})$ will be equal to $\qdisc(\cl{P}_{\cl{S}})$ for some set system $\cl{S}$. However, this is not the case for a general projection system which may include non-commuting elements. Hence, the complexity of the notion of quantum discrepancy also lies in the non-commutativity inside $\cl P$.

Now, let us justify the term ``quantum'' in ``quantum discrepancy''. In classical physics, the phase space of a system is described by a set which is finite if the system has finite states.  An observable of this system is a real-valued function on the phase space. In quantum physics, the phase space is replaced by a Hilbert space. If the classical phase space has $N$ elements, the quantum phase space would be $\bb C^N$. The quantum counterpart of an observable $f$ is a Hermitian operator with the same spectrum as the image of $f$ (see \cite{WeaverBook} for more explanation). Now, note that in the current quantization process of combinatorial discrepancy, subsets of $[N]$ are replaced by orthogonal projections which are in one-to-one correspondence with subspaces of $\bb C^N$, and functions of the form \linebreak[4] $\rchi: [N]\to\{\pm 1\}$ are replaced by Hermitian matrices with eigenvalues in $\{\pm 1\}$. 

In the preface of \cite{WeaverBook} Weaver refers to the two above-mentioned ideas as general principles of quantization. It is worthwhile to see this in his own words:  
\begin{quotation}
... We have now reached a point where it is possible to give a simple, unified approach to the general concept of quantization in mathematics. 

 The fundamental idea of mathematical quantization is that sets \emph{are replaced by Hilbert spaces}. Thus, we regard lattice operations (join, meet, orthocomplement) on subspaces of a Hilbert space as corresponding to set-theoretic operations (union, intersection, complement) on subsets of a set. This already allows one to determine quantum analogs of some simple structures. But the real breakthrough is the fact that the quantum version of a complex-valued function on a set is an operator on a Hilbert space. ...

With more work the analogy can be pushed even further. At each step one must formulate the given classical notion in just the right way
to obtain a viable quantum version. This sometimes requires significant creativity. However, as it is done in case after case, general quantization principles emerge.
\end{quotation}
 
The remaining parts of this paper are organized as follows. The next subsection is dedicated to the statement of the main results. In order to make the current subsection more coherent, discussions about equivalent forms of quantum discrepancy are postponed to Section~\ref{Sec: EquivForms}. We prove the main results in Sections~\ref{Sec: GeneralProofs}, and~\ref{Sec: SetSysProofs}. Finally, Section~\ref{Sec: ConcRems} is devoted to some concluding remarks.  
\subsection{Main Results}
The main results of this paper can be divided into two categories: quantum discrepancy of a general projection system, and quantum discrepancy of a set system. Note that, to be more precise, we have to use a phrase like ``quantum discrepancy of a projection system which is associated with a set system'' instead of  ``quantum discrepancy of a set system''. But, for more convenience, we use the latter.

 In the case of a general projection system, we first investigate upper bounds. For a set system $\cl S$ of size $M$ in $[N]$, a trivial upper bound for $\disc(\cl S)$ is $N$. We prove the same bound holds for quantum discrepancy.
 \begin{theorem}[Trivial Upper Bound]
 	\label{Thm: TrivBound}
 	For each projection system $\cl{P}$ in $\bb{C}^N$,
 	\begin{align*}
 	\qdisc(\cl{P})\leq N.
 	\end{align*}
 \end{theorem}
 Another general upper bound of $\disc(\cl S)$ can be obtained by coloring the points of $[N]$, randomly. It can be shown that 
 if $\cl{S}$ is a set system consisting of $M$ subsets of $[N]$, $\disc(\cl{S}) = O(\sqrt{N\log M})$. In particular for $M= O(N)$, $\disc(\cl{S})=O(\sqrt{N\log N})$. 
 	
 We recall that for two functions $f$ and $g$ on $\bb N\times\bb N$,  $f(m,n)=O(g(m,n))$ means that there exists a constant $c>0$ independent of $m$ and $n$, such that $|f(m,n)|\leq c|g(m,n)|$ for all but a finite set of values of $(m,n)\in\bb N\times\bb N$.
 
  Using random quantum colorings, we obtain an upper bound for $\qdisc(\cl P_{\cl S})$ in terms of $N$ and $M$.
\begin{theorem}
	\label{Thm: UpperBound}
	Suppose $\cl{P}$ is an $M$\!-element projection system in $\bb{C}^N$. Then
	\begin{align*}
	\qdisc(\cl{P}) = O(\sqrt{N+\log M}).
	\end{align*}
\end{theorem}
 The random coloring bound of $\disc(\cl S)$, i.e. $O(\sqrt{N\log M})$,  is not tight. Spencer, in \cite{Spencer}, obtained another upper bound which is optimal. He proved that when $M\geq N$, 
 	\begin{align*}
 	\disc(\cl{S}) = O\left(\sqrt{N \log\left(2M/N\right)}\right).
 	\end{align*}
 	In particular, if $M = O(N)$, $\disc(\cl{S})$ would be of order $O(\sqrt{N})$. A probabilistic proof for the optimality of Spencer's bound is sketched in~\cite{Matousek} (Exercise 1, Section 4.1). It is based on a random set system $\cl{S}$ and is designed to prove that if $c>0$ is a sufficiently small constant  and $N\leq M\leq 2^{cN}$, there exists a set system $\cl{S}$ on $[N]$ with $M$ elements such that $\disc(\cl{S})\geq c\sqrt{N\log({2M}/{N})}$. 

   In case of a projection system, its size is not bounded in terms of the dimension of the ground space. Hence, by Theorem~\ref{Thm: TrivBound}, the bound given in Theorem~\ref{Thm: UpperBound} can not be tight for all values of $N$ and $M$. However, by choosing a random projection system with an appropriate distribution, we prove that this bound is tight in the worst case, provided that $M$ is neither too large nor too small with respect to $N$.
\begin{theorem}
	\label{Thm: Tightness}
	Suppose	$M$ is neither too small nor too large with respect to $N$, such that	
	\begin{align*}
	aN\leq M\quad,\quad \log M\leq bN,
	\end{align*}
	for appropriate constants $a,b>0$, which are independent of $M$ and $N$. There exist an $M$-element projection system $\cl{P}$ in $\bb{C}^N$ and a constant $c>0$, which depends only on $b$, with the property that
	\begin{align*}
	\qdisc(\cl{P})\geq c\sqrt{N+\log M}.
	\end{align*}
\end{theorem}
In the procedure of proving Theorem~\ref{Thm: UpperBound} it becomes clear that this partially tight bound occurs with probability at least $1/2$. Hence, Theorem~\ref{Thm: UpperBound} is of constructive importance. If we generate multiple samples of the specified random quantum coloring, we have, with high probability, a quantum coloring that satisfies the given bound. 

After this general consideration, we focus on a projection systems $\cl P_{\cl S}$ corresponding to a set system $\cl S$ of size $M$ in $[N]$. For this special kind of projection systems both $\disc(\cl{S})$ and  $\qdisc(\cl{P}_{\cl{S}})$ are defined.  Since each combinatorial coloring is a quantum coloring, 
\begin{align}
\label{rel: TrivDiscQuant}
\qdisc(\cl P_{\cl S})\leq\disc(\cl S).
\end{align}
We prove two other relations between these two quantities. First, we investigate the case where $\qdisc(\cl P_{\cl S})=0$, and obtain the following result.
\begin{theorem}
	\label{Thm: ZerosDiscQDisc}
	Suppose $\cl{S}$ is a set system on $[N]$, and $\cl{P}_{\cl{S}}$ is its corresponding projection system. $\qdisc(\cl{P}_{\cl{S}})=0$ if and only if $\disc(\cl S)=0$.
\end{theorem}
For the remaining set systems we bound $\qdisc(\cl P_{\cl S})$ and $\disc(\cl S)$ in terms of each other in a way different from~\eqref{rel: TrivDiscQuant}.
\begin{theorem}
	\label{Thm: QuantCombThm2}
	Suppose $\cl{S}$ is a set system with $M$ subsets of  $[N]$, and $\cl{P}_{\cl{S}}$ is its corresponding projection system. If  $\qdisc(\cl{P}_{\cl{S}})\neq 0$, then for some constant $c>0$, independent of $M$, $N$, and $\cl S$, it holds that
	\begin{align*}
	\disc(\cl{S})<\left(2c\max\left[\sqrt{\log(2M)}\,\cc\,\frac{\log(2M)}{\qdisc(\cl{P}_{\cl{S}})}\right]+1\right)\qdisc(\cl{P}_{\cl{S}}).
	\end{align*}
\end{theorem}

This theorem shows that $\disc(\cl S)$ which is not smaller than $\qdisc(\cl P_{\cl S})$ can not be arbitrarily larger. This illustrates an important aspect of the study of quantum discrepancy to combinatorial discrepancy theory, i.e. an upper bound for $\qdisc(\cl P_{\cl S})$ leads to an upper bound for $\disc(\cl S)$. The idea of the proof is that a quantum coloring with low quantum discrepancy gives us a low-discrepancy coloring. The next example demonstrates another application of Theorem~\ref{Thm: QuantCombThm2}.
\begin{example}
	Suppose $\cl{S}$ is the set system of all arithmetic progressions in $[N]$. As mentioned before, $\disc(\cl S)=\Omega(N^{1/4})$. Thus, 
	$\qdisc(\cl P_{\cl S})\neq 0$. It can be shown that in this special case $M= O(N^3)$, so $\log (2M) = O(\log N)$. From Theorem~\ref{Thm: QuantCombThm2} we conclude that for some $c'>0$, we have $(c'\sqrt{\log N}+1)\qdisc(\cl P_{\cl S})=\Omega(N^{1/4})$. Hence,
	\begin{align*}
	\qdisc(\cl{P}_{\cl{S}})=\Omega(N^{\frac{1}{4}-\epsilon}),\quad \forall \epsilon>0\cdot
	\end{align*}
	And since $\disc(\cl S)=O(N^{1/4})$ 	(see~\cite{MatousekSpencer}), it follows that for each $\epsilon>0$, 
	\begin{align*}
	c_1N^{\frac{1}{4}-\epsilon}\leq\qdisc(\cl{P}_{\cl{S}})\leq c_2N^{\frac{1}{4}},
	\end{align*}
	where $c_1=c_1(\epsilon)$ and $c_2$ are suitable positive constants.
\end{example}

\section{Equivalent Forms of Quantum Discrepancy}
\label{Sec: EquivForms}
	

\subsection{A Probabilistic Approach to Combinatorial Discrepancy}
\label{Sub: ProbApp}
A probabilistic interpretation of quantum discrepancy when restricted to set systems is explained in this subsection. The core idea behind the scenes is to color points of the ground set $[N]$ randomly according to the distribution of a determinantal process. 

Note that a coloring of $[N]$ is uniquely determined by a subset of $[N]$, that is the set of, say, its red points. Hence, a random coloring of $[N]$ is, in fact, a random subset of this set. Such a random object, i.e a random subset of $[N]$, is called a \emph{simple point process} on $[N]$ (the word ``simple'' indicates that the multiplicity of each point of $[N]$ in any realization is at most one).

\begin{definition}
	\label{Def: DetPross}
	A simple point process $\mathfrak{X}$ on $[N]$ is said to be a \emph{determinantal process}, if there is a matrix $K\in\cl{M}_N(\bb C)$ such that for every $\{i_1,i_2,\ldots,i_k\}\subseteq[N]$
	\begin{align}
	\label{eq: DetProssDef}
	\bb{P}\left[i_1,i_2,\ldots,i_k\in\mathfrak{X}\right]=\det\left[\left(K_{ij}\right)_{i,j\in\{i_1,i_2,\ldots,i_k\}}\right]\cdot
	\end{align}
   $K$ is called the \emph{kernel} of  $\mathfrak{X}$.
\end{definition}                                       
The following proposition provides a necessary and sufficient condition for a Hermitian matrix $K$ to be the kernel of a determinantal process (see e.g. Theorem 4.5.5 in~\cite{GAF}).
\begin{proposition}
	\label{Prop: HermitianKernel}
	Suppose $K$ is a Hermitian matrix. $K$ is the kernel of a determinantal process if and only if its spectrum is a subset of $[0,1]$.
\end{proposition} 

A determinantal process $\mathfrak{X}$ with a Hermitian kernel $K$ on $[N]$ has negative correlations since for any distinct values of $i,j \in [N]$
 \begin{align*}
\bb{P}[i,j\in\mathfrak{X}] &= \bb P[i\in\mathfrak{X}]\bb P[j\in\mathfrak{X}]-|K_{ij}|^2\\
&\leq\bb P[i\in\mathfrak{X}]\bb P[j\in\mathfrak{X}]\cdot
\end{align*}
More generally, using Koteljanski\u\i's inequality (see e.g.~\cite{Horn}, Theorem 7.8.9), it can be proved that any two disjoint sets $I,J\subseteq[N]$ repel each other in the sense that \begin{align*}
\bb{P}[I\cup J\subseteq\mathfrak{X}]\leq\bb{P}[I\subseteq\mathfrak{X}]\bb{P}[J\subseteq\mathfrak{X}]\cdot
\end{align*}
Therefore, it is expected that determining red points of a coloring according to the law of a determinantal process would prevent them from accumulating in one or more members of a given set system. 

Now, suppose $\rchi$ is a quantum coloring.  Then, $K\coloneqq(\rchi+I)/2$ is a Hermitian matrix with eigenvalues in $\{0, 1\}$. In fact, $K$ is an orthogonal projection. On the other hand, for an arbitrary orthogonal projection $K$, the matrix $2K-I$ is a quantum coloring. Hence, by Proposition~\ref{Prop: HermitianKernel}, there is a one-to-one correspondence between the set of $N\times N$ quantum colorings and the set of \emph{determinantal projection processes} on $[N]$, i.e. determinantal processes that their kernels are orthogonal projections. The mentioned probabilistic interpretation of quantum discrepancy of a set system is based on this correspondence.

\begin{lemma}
	\label{lem: ExpToMat}
	Suppose $\cl P_{\cl S}$ is the projection system corresponding to a set system $\cl S$ in $[N]$. Then,
	\begin{align*}
	\qdisc(\cl P_{\cl S}) = \min_{\mathfrak{X}}\max_{S\in\cl S}	\left[\bb{E}\left[\left({2\mathfrak{X}}(S)-|S|\right)^2\right]\right]^\frac{1}{2},
	\end{align*} 
while the minimum is taken over all determinantal projection processes on $[N]$.
\end{lemma}
\begin{proof}
	According to the explanation given before the lemma, it is enough to prove that for each $S\subseteq[N]$, 
	\begin{align}
	\label{eq: MatToExp}
	\tr^2(\rchi P_S)+\tr\left(\rchi[\rchi,P_S]P_S\right) =\bb{E}\left[\left({2\mathfrak{X}}(S)-|S|\right)^2\right]\cc
	\end{align}
	where $\rchi$ is an arbitrary quantum coloring and $\mathfrak{X}$ is the determinantal projection process with kernel $K = (\rchi+I)/2$. We start from the right-hand side. By definition of  $P_S$, we have $|S| = \tr(P_S)$. Thus,
	\begin{align}
	\label{eq: ExpToMat1}
	\bb E\left[\left(2\mathfrak X(S)- |S|\right)^2\right]&=4\bb E\left[\mathfrak X^2(S)\right]-4\tr(P_S)\bb E\left[\mathfrak X(S)\right]+\tr^2(P_S)\cdot
	\end{align}
	We compute $\bb E{[\mathfrak{X}(S)]}$ and $\bb E\left[\mathfrak{X}^2(S)\right]$. Suppose that the $i$-th diagonal element of $P_S$ is $p_{ii}$. If $\textbf{1}_E$ represents the indicator function of an arbitrary event $E$, then for any $S\subseteq[N]$ it holds that
	\begin{align}
	\label{eq: SumOfIndic}
	\mathfrak{X}(S) = \sum_{i\in S}\textbf{1}_{\{i\in\mathfrak X\}}.
	\end{align}
	Hence,
	\begin{align}
	\bb{E}\left[\mathfrak X(S)\right] &=\sum_{i\in S}\bb{P}\left[i\in\mathfrak X\right]\tag{by \eqref{eq: SumOfIndic}, linearity of expectation}\nonumber\\
	&=\sum_{i=1}^{ N}\bb P\left[i\in\mathfrak X\right]\textbf{1}_{\{i\in S\}}\nonumber.
	\end{align}
	According to the definition of $P_S$, we have $\textbf{1}_{\{i\in S\}}=p_{ii}$. Thus,
	\begin{align}
	\bb{E}\left[\mathfrak X(S)\right]&= \sum_{i=1}^{ N}K_{ii}p_{ii}\nonumber\tag{by \eqref{eq: DetProssDef}}\\
	&=\tr(KP_S)\cdot\label{eq: ExpToMat2}
	\end{align}
	Using~\eqref{eq: SumOfIndic} and linearity of the expectation once again, we obtain
	\begin{align}
	\bb{E}\left[\mathfrak{X}^2(S)\right]
	&=\sum_{i\in S}\bb{P}\left[i\in\mathfrak X\right]+\sum_{\substack{S\times S\\i\neq j}}\bb{P}\left[i,j\in\mathfrak X\right]\nonumber\\
	& = \tr(KP_S)+\sum_{\substack{S\times S\\i\neq j}}\bb{P}\left[i,j\in\mathfrak X\right]\cdot\label{eq: ExpToMat3}
	\end{align}
	On the other hand,
	\begin{align}
	\sum_{\substack{S\times S\\i\neq j}}\bb P\left[i, j\in\mathfrak X\right]&=\sum_{1\leq i\neq j\leq N}\bb P\left[i, j\in\mathfrak X\right]\textbf{1}_{\{i \in S\}}\textbf{1}_{\{j\in S\}}\nonumber\\
	&=\sum_{1\leq i\neq j\leq N}\left[K_{ii}K_{jj}-K_{ij}K_{ji}\right]p_{ii}p_{jj}\nonumber\tag{by \eqref{eq: DetProssDef}}\\
	& = \sum_{1\leq i\neq j\leq N}K_{ii}p_{ii}K_{jj}p_{jj}-\sum_{1\leq i\neq j\leq N}|K_{ij}|^2p_{ii}p_{jj}\nonumber\tag{$K$ is Hermitian}\\
	& = \tr^2(KP_S)-\sum_{1\leq i\leq N}K_{ii}^2p_{ii}-\sum_{1\leq i\neq j\leq N}|K_{ij}|^2p_{ii}p_{jj}\nonumber\\
	&=\tr^2(KP_S)-\sum_{1\leq i,j\leq N}|K_{ij}|^2p_{ii}p_{jj}\nonumber\\
	&=\tr^2(KP_S)-\tr\left((KP_S)^2\right)\cdot\label{eq: ExpToMat4}
	\end{align}    
	Combining \eqref{eq: ExpToMat1}, \eqref{eq: ExpToMat2}, \eqref{eq: ExpToMat3}, and \eqref{eq: ExpToMat4} results in the following equation.
	\begin{align}
	\label{eq: ExpToKer}
      \bb E\left[(\left(2\mathfrak X(S)- |S|\right)^2\right]=\tr^2(2KP_S- P_S)+4\tr\left(KP_S(I-KP_S)\right).
	\end{align}
	 It is straightforward to check that for $K = (\rchi+I)/2$, and each $S$,
	\begin{align*}
	\tr^2(2KP_S- P_S)+4\tr\left(KP_S(I-KP_S)\right) = \tr^2(\rchi P_S)+\tr\left(\rchi[\rchi,P_S]P_S\right),
	\end{align*}
	which together with \eqref{eq: ExpToKer} gives \eqref{eq: MatToExp}.
	
\end{proof}

\begin{remark}
	By the bias-variance decomposition
	\begin{align*}
	\bb{E}\left[\left(2\mathfrak{X}(S)-|S|\right)^2\right]&=\big[\bb{E}\left[2\mathfrak{X}(S)\right]-|S|\big]^2+\rm{Var}\left[2\mathfrak{X}(S)\right].
	\end{align*}
	If $\mathfrak X$ is the uniform independent coloring of points, and so has kernel $I/2$ (which is not a projection), the first summand on the right-hand side will be zero. Moreover, for a deterministic coloring $\mathfrak{X}$ the second summand vanishes. Lemma~\ref{lem: ExpToMat} indicates that the quantum discrepancy of a set system is small, if there exists a determinantal projection process on the ground set that makes both summands small enough for each member of that set system.  \\	
	Note that the set of orthogonal projections in $\bb C^N$ is the set of extreme points of the convex set of positive semi-definite elements of $\cl{M}_N(\bb C)$ with eigenvalues in $[0,1]$. In other words, according to Proposition~\ref{Prop: HermitianKernel}, kernels of all determinantal projection processes constitute the extreme points of the set of kernels of all determinantal processes with a Hermitian kernel. This makes a similarity between the set of (combinatorial) colorings and quantum colorings. Each coloring of $[N]$ can be assigned uniquely to an $N$-tuple of $0$s and $1$s which determines its distribution: the $i$-th coordinate shows the probability of coloring the element $i$ red. These tuples are extreme points of the unit cube $[0,1]^N$ which can be viewed as the set of laws of all methods for coloring points of $[N]$ independently: in $(p_1,p_2,\ldots,p_N)\in [0,1]^N$, $p_i$ is the chance of $i$ to be colored red. Moreover, just as the number of red points in each coloring is fixed, it can be shown that determinantal projection processes are exactly those determinantal processes with Hermitian kernels that have a fixed number of points (see the second part of Proposition~\ref{prop: Bernoullis} which is a part of Theorem 4.5.3 in~\cite{GAF}).
\end{remark}

\subsection{A Measure of Non-Commutativity}
\label{SubSec: NormCommutForm}
It is explained in Subsection~\ref{Sub: Definition of QuantDisc} that the term $\tr\left(\rchi[\rchi,P]P\right)$ encodes the non-commutative nature of quantum discrepancy. The following lemma provides an algebraically simple relation between $\tr\left(\rchi[\rchi,P]P\right)$ and the commutator of $\rchi$ and $P$, which illustrates the existing non-commutativity more explicitly.

\begin{lemma}
	\label{lem: NormTrace}
	Suppose $P$ and $\rchi$ are an orthogonal projection and a quantum coloring, respectively. Then,
	\begin{align*}
	\tr\left(\rchi[\rchi,P]P\right)=\frac{1}{2}\|[\rchi,P]\|_2^2\cdot
	\end{align*}
\end{lemma}
\begin{proof}
	The following calculation results in what we want.
	\begin{align*}
	\|[\rchi,P]\|_2^2 & = \tr([\rchi,P]^*[\rchi,P])\\
	& =\tr\left((P\rchi-\rchi P)(\rchi P-P\rchi))\right)\\
	&=\tr\left(P\rchi^2 P-(P\rchi)^2-(\rchi P)^2+\rchi P^2\rchi\right)\\
	&=2\tr\left(\rchi^2P^2-(\rchi P)^2\right)\tag{trace is cyclic}\\
	&=2\tr\left(\rchi(\rchi P-P\rchi)P\right)\\
	&=2\tr(\rchi[\rchi,P]P)\cdot
	\end{align*}
\end{proof}

\begin{corollary}
	\label{Col: NormFormulation}
	Suppose $\cl P$ is a projection system in $\bb C^N$. Then,
	\begin{align*}
	\qdisc(\cl P)=\min_{\rchi}\max_{P\in\cl P}\left[\tr^2(\rchi P)+\frac{1}{2}\|[\rchi,P]\|_2^2\right]^{1/2},
	\end{align*}
	where the minimum is taken over the set of quantum colorings.
\end{corollary}
Now, quantum discrepancy can be interpreted more simply. As mentioned in Subsection~\ref{Sub: Definition of QuantDisc}, for a diagonal projection $P_S$ and a diagonal coloring  $\rchi$, $\tr(\rchi P_S) = \rchi(S)$. Therefore, the term $|\tr(\rchi P)|$ is common between the two kinds of discrepancies. The additional term, $\frac{1}{2}\norm{[\rchi,P]}{2}^2$, appears as a measure of the non-commutativity arisen as the result of replacing a (diagonal) coloring and a set with a quantum coloring and an orthogonal projection, respectively.

\begin{remark}
	The specific combination of $|\tr(\rchi P)|$ and $\norm{[\rchi,P]}{2}$ used to define quantum discrepancy, i.e. $\left[\tr^2(\rchi P)+\frac{1}{2}\norm{[\rchi,P]}{2}^2\right]^{1/2}$, has the advantage of having a probabilistic interpretation as explained in the previous subsection. However, since 
	\begin{align*}
	\frac{a+b}{\sqrt{3}}\leq\sqrt{a^2+\tfrac{1}{2}\,b^2}\leq a+b,\quad \text{for all $a,b\geq 0$},
	\end{align*}
	$\qdisc(\cl P)$ is essentially the same quantity as the simpler and possibly more intuitive formulation $\min_{\rchi}\max_{P\in\cl P}\big[|\tr(\rchi P)|+\norm{[\rchi,P]}{2}\big]$. 
\end{remark}

\section{General Upper Bounds}
\label{Sec: GeneralProofs}
For a set $S\subseteq [N]$ and a determinantal process $\mathfrak{X}$, we have $\big|2\mathfrak{X}(S)-|S|\big|\leq N$. Hence, $\qdisc(\cl{P}_{\cl{S}})\leq N$ for each set system $\cl S$. Theorem~\ref{Thm: TrivBound} states that this bound also holds for an arbitrary projection system. 
\begin{proof}[Proof of Theorem~\ref{Thm: TrivBound}]
	We show that if $\rchi,P\in\cl{M}_N(\bb C)$ are a quantum coloring and an orthogonal projection, respectively, then 
	   \begin{align*}
	\tr^2\left(\rchi P\right)+\tr\left(\rchi[\rchi,P]P\right)\leq N^2,
	\end{align*}  
		and this results in what we want. According to the properties of $\rchi$ and $P$, 
		\begin{align}
		\label{eq: CommutSubtrac}
		 \tr\left(\rchi[\rchi,P]P\right)=\tr\left(P-(\rchi P)^2\right).
		 \end{align}
		  Therefore, we show that
        \begin{align*}
	\tr^2\left(\rchi P\right)+\tr\left(P-(\rchi P)^2\right)\leq N^2.
	\end{align*}  
 Because of the cyclic property of the trace function, it is enough to prove this inequality for a diagonal $\rchi$. To be more precise, consider a spectral decomposition of $\rchi$, i.e. a representation such as $\rchi = UDU^*$, where $D$ is diagonal and $U$ is unitary. Then, for instance, $\tr(\rchi P) = \tr\left(D(U^*PU)\right)$, and $U^*PU$ is an orthogonal projection.\\
 By computing diagonal entries of $\rchi P$ and $(\rchi P)^2$ it can be shown that
	\begin{align*}
	\tr^2\left(\rchi P\right)-\tr\left((\rchi P)^2\right)=2\sum_{1\leq i<j\leq N}\left[\rchi_{ii}\rchi_{jj}\left(P_{ii}P_{jj}-|P_{ij}|^2\right)\right].
	\end{align*}   
 The orthogonal projection $P$ is positive semi-definite and its operator norm is equal to $1$.  Hence, $P_{ii}P_{jj}-|P_{ij}|^2$ is in $[0,1]$. Diagonal entries of $\rchi$ are $-1$ or $1$, and we conclude that
	\begin{align*}
	\tr^2\left(\rchi P\right)-\tr\left((\rchi P)^2\right)&\leq 2\sum_{1\leq i<j\leq N}\!1\\
	&=N^2-N.
	\end{align*}
	At last, $\tr^2\left(\rchi P\right)+\tr\left(P-(\rchi P)^2\right)\leq N^2-N+N=N^2$.
	
\end{proof}

\begin{remark}
	With a similar method, the bound given in Theorem~\ref{Thm: TrivBound} can be improved. In fact, it holds that
	\begin{align*}
		\qdisc(\cl P)\leq\max_{P\in\cl P}\ra(P),
	\end{align*}
	which results in $	\qdisc(\cl P)\leq N$.
\end{remark}

\subsection{Proof of Theorem~\ref{Thm: UpperBound}}
Instead of proving Theorem~\ref{Thm: UpperBound} directly, we prove Theorem~\ref{thm: quantrandcol} which results in the former. Our approach to give an upper bound for the quantum discrepancy of a general projection system $\cl{P}$ in $\bb{C}^N$ is to investigate the random variable $\left[\tr^2(\rchi P)+\tr\left(\rchi[\rchi,P]P\right)\right]^{1/2}$, where $\rchi$ is a  random quantum coloring with a certain distribution (Definition~\ref{dfn: quantrandcol}) and $P\in\cl{P}$. For more convenience we name this random variable $\rchi(P)$. As mentioned before,
\begin{align} 
\label{eq: ChiP}
\rchi(P)=\left[\tr^2(\rchi P)+\tr\left(P-(\rchi P)^2\right)\right]^{\frac{1}{2}}.
\end{align}
\begin{definition}
	\label{dfn: quantrandcol}
	By a \emph{random balanced quantum coloring} we mean a random matrix \linebreak $\rchi= UDU^*$,  where $U$ has Haar distribution on the group of $N\times N$ unitary matrices, $\bb{U}(N)$, and
	\begin{align*}
	D\coloneqq\left[\begin{array}{c|c}
	I_{\lfloor\frac{N}{2}\rfloor}&\mathbf{0}\\
	\hline
	\mathbf{0}&-I_{\lceil\frac{N}{2}\rceil}\end{array}\right]\cdot
	\end{align*}
\end{definition}
 
\begin{theorem}[Quantum Random Coloring]
	\label{thm: quantrandcol}
	Suppose $\cl{P}$ is an $M$-element projection system in $\bb{C}^N$. There exists a constant $c>0$, independent of $N$, $M$ and $\cl P$, such that for
	\begin{align*}
	\Delta_P\coloneqq \sqrt{2}\left[\sqrt{\frac{1}{c}\log(8M)+\frac{N^2}{N^2-1}\ra(P)- \frac{N}{N^2-1}\ra^2(P)}+\frac{\ra(P)}{N}\right]\cc
	\end{align*}
	and for large enough values of  $N$, the probability that a random balanced quantum coloring $\rchi$ satisfies simultaneously all inequalities
	\begin{align*}
	\rchi(P)\leq \Delta_P\cc\quad P\in\mathcal P
	\end{align*}
	is at least $1/2$. In particular, $\qdisc(\cl{P})=O(\sqrt{N+\log M})$.
\end{theorem}

Proof of this theorem is mainly based on the concentration of $\tr(\rchi P)$ and $\tr\!\left(P\!-\!(\rchi P)^2\right)$ around their means. We will use the following proposition which is the same as Corollary 4.4.31 in \cite{AndGuiZei} but in our words.
\pagebreak[4] 
\begin{proposition}
	\label{prop: ConcProp}
	Assume we are given deterministic matrices $X_1,\ldots,X_k\in \cl{M}_N(\mathbb C)$ and a constant $\sigma$ that controls all singular values of these matrices from above. Let \linebreak $p\coloneqq p(x_1,x_2,\ldots,x_{k+2})$ be a polynomial of $k+2$ non-commutative variables with complex coefficients. For $X\in \bb{U}(N)$ define $f(X)=\tr\big( p(X,X^*,X_1,X_2,\ldots,X_k)\big)$. Then, there are positive constants $N_0=N_0(p)$ and $c=c(p,\sigma)$ such that for any	$\delta>0$ and $N>N_0$,
	\begin{align*}
	\bb{P}\Big[\big|f(X)-\bb{E}[f(X)]\big|\geq\delta N\Big]\leq 2\exp\left(-cN^2\delta^2\right)\cc
	\end{align*}
	in which, $\bb{P}$ is the unique Haar probability measure on $\bb{U}(N)$ and $\bb{E}[\cdot]$ is the expected value with respect to $\bb{P}$.
\end{proposition}

Another component of the proof is estimating $\bb{E}\left[\tr(\rchi P)\right]$ and $\bb{E}\left[\tr(P-(\rchi P)^2)\right]$. 
\begin{lemma}
	If  $\rchi$ is a random balanced quantum coloring, and $P$ is an orthogonal projection in $\mathbb C^N$, then	
	\begin{align}
	\label{chiP}
	&\Big|\bb{E}\big[\tr(\rchi P)\big]\Big|\leq\frac{\ra(P)}{N}\cc\\
	\label{chiP2}
	&\bb{E}\Big[\tr\big(P-(\rchi P)^2\big)\Big]\leq \frac{N^2}{N^2-1}\ra(P)- \frac{N}{N^2-1}\ra^2(P)\cdot
	\end{align}
\end{lemma}
\begin{proof}
	We compute exact values of $\bb{E}\big[\tr(\rchi P)\big]$ and $\bb{E}\left[\tr\big(P-(\rchi P)^2\big)\right]$. According to the distribution of $\rchi$, we assume, without loss of generality, that $P$ is diagonal. More precisely, if $V$ and $\Pi$ are respectively a unitary and a diagonal matrix such that $P=V\Pi V^*$ is a spectral decomposition of $P$, then by the cyclic property of the trace function
	\begin{align*}
	\tr(\rchi P) = \tr(UDU^*V\Pi V^*) =\tr\left((V^*U)D(V^*U)^*\Pi\right)\cdot
	\end{align*}
	Haar measure is invariant under left multiplication by a unitary matrix, so $(V^*U)D(V^*U)^*$ has the same distribution as $\rchi$.  
	
	\vspace*{.4cm} \noindent
	$\boldsymbol{\bb{E}\big[\tr(\rchi P)\big]:}$
	\begin{align*}
	(\rchi P)_{ii} = \sum_{j=1}^NU_{\!ij}D_{\!jj}\overline{U}_{\!ij}P_{\!ii}\tag{$D, P$ are diagonal}\cdot
	\end{align*}
	Thus, $\bb{E}\left[(\rchi P)_{ii}\right]=\sum_{j=1}^N\bb{E}\left[\left|U_{\!ij}\right|^2\right]D_{\!jj}P_{\!ii}$. Since Haar measure is invariant under left and right multiplication by a unitary matrix, all entries of $U$ are identically distributed. As a result, 
	$\bb{E}\left[\left|U_{\!ij}\right|^2\right]=1/N$. Therefore,
	\begin{align*}
	\bb{E}\big[\tr(\rchi P)\big]=\sum_{i,j}\frac{1}{N}D_{jj}P_{ii}
	=\frac{1}{N}\tr(D)\tr(P)=\left\{\begin{array}{lr}
	0\cc&\text{even $N$}\\
	-\dfrac{\ra(P)}{N}\cc&\text{odd $N$}
	\end{array}\cdot
	\right.
	\end{align*}
	Now, \eqref{chiP} is clear.
	
	\pagebreak
	\vspace{.4cm}\noindent
	$\boldsymbol{\bb{E}\Big[\tr\big(P-(\rchi P)^2\big)\Big]:}$\vspace{.1cm}
	
	\noindent
	According to Lemma~\ref{lem: NormTrace}, and \eqref{eq: CommutSubtrac}, it holds that $\tr\big(P-(\rchi P)^2\big)=\frac{1}{2}\norm{[\rchi,P]}{2}^2$. Therefore,
	$\tr\big(P-(\rchi P)^2\big)=\frac{1}{2}\sum_{i,j}\left|([\rchi,P])_{ij}\right|^2$. We have
	\begin{align*}
		([\rchi,P])_{ij}=\rchi_{ij}\left(P_{jj}-P_{ii}\right)=\left\{\begin{array}{ll}
		\rchi_{ij}\cc\quad&P_{ii}=0, P_{jj}=1\\
		-\rchi_{ij}\cc\quad&P_{ii}=1, P_{jj}=0\\
		0,&\text{otherwise}
		\end{array}\right.\cdot
	\end{align*}
	Therefore, $\norm{[\rchi,P]}{2}^2=2\sum_{i,j:P_{ii}=1, P_{jj}=0}\left|\rchi_{ij}\right|^2$. We compute $\left|\rchi_{ij}\right|^2$.
	\begin{align*}
	\left|\rchi_{ij}\right|^2=\rchi_{ij}\cdot\overline{\rchi}_{\!ij}&=\sum_{k=1}^NU_{ik}D_{kk}\overline{U}_{\!jk}\cdot\sum_{l=1}^N
	\overline{U}_{\!il}D_{ll}U_{jl}\tag{$D$ is diagonal}\\
	&=\sum_{k,l}U_{ik}U_{jl}\overline{U}_{\!il}\overline{U}_{\!jk}D_{kk}D_{ll}.
	\end{align*} 
	 Using Proposition 4.2.3 in~\cite{HiaiPetz}, we know that for $i\neq j$:
	\begin{align*}
	\bb E\big[U_{ik}U_{jl}\overline{U}_{\!il}\overline{U}_{\!jk}\big]&=\left\{\begin{array}{ll}
      \bb{E}\left[\left|U_{ik}\right|^2\left|U_{jk}\right|^2\right]\cc\qquad& k=l\\
      \\
      \bb{E}\left[U_{ik}U_{jl}\overline{U}_{\!il}\overline{U}_{\!jk}\right]\cc\qquad&k\neq l
	\end{array}\right.\\
	&=\left\{\begin{array}{ll}
	\dfrac{1}{N(N+1)}\cc\qquad& k=l\\
	\\
	-\dfrac{1}{N(N^2-1)}\cc\qquad&k\neq l
	\end{array}\right..
	\end{align*}
	Therefore,
	\begin{align*}
		\bb E\big[\left|\rchi_{ij}\right|^2\big]&=\sum_{k=1}^N\dfrac{1}{N(N+1)}D_{kk}^2-\sum_{k\neq l}\dfrac{1}{N(N^2-1)}D_{kk}D_{ll}\\
		&=\sum_{k=1}^N\dfrac{1}{N(N+1)}-\dfrac{1}{N(N^2-1)}\sum_{k\neq l}D_{kk}D_{ll}\tag{$D_{kk}=\pm1$}\\
		&=\dfrac{1}{N+1}-\dfrac{1}{N(N^2-1)}\left(\tr^2(D)-N\right)\\
		&=\left\{
			\begin{array}{ll}
			\dfrac{N}{N^2-1}\cc\qquad&\text{even $N$}\\
			\\
			\dfrac{1}{N}\cc\qquad&\text{odd $N$}
			\end{array}
			\right..
	\end{align*}
We obtain
	\begin{align*}
	\bb{E}\big[\tr\left(P-(\rchi P)^2\right)\big]&=\frac{1}{2}\,\bb E\big[\norm{[\rchi,P]}{2}^2\big]\\
	&=\sum_{i,j:P_{ii}=1,P_{jj}=0}\bb E\left[\left|\rchi_{ij}\right|^2\right]\\
	&=\left\{\begin{array}{ll}
	\dfrac{N}{N^2-1}\cdot\sum_{i,j:P_{ii}=1,P_{jj}=0}1\cc\qquad&\text{even $N$}\\
	\\
	\dfrac{1}{N}\cdot\sum_{i,j:P_{ii}=1,P_{jj}=0}1\cc\qquad&\text{odd $N$}
	\end{array}
	\right.\\
	&=\left\{\begin{array}{ll}
	\dfrac{N}{N^2-1}\cdot\ra(P)\!\left(N-\ra(P)\right)\cc\qquad&\text{even $N$}\\
	\\
	\dfrac{1}{N}\cdot\ra(P)\!\left(N-\ra(P)\right)\cc\qquad&\text{odd $N$}
	\end{array}
	\right.\\
	&\leq\dfrac{N}{N^2-1}\cdot\ra(P)\!\left(N-\ra(P)\right)\tag{$N\geq 2$}\\
	&=\frac{N^2}{N^2-1}\ra(P)- \frac{N}{N^2-1}\ra^2(P)\cdot
	\end{align*}
	Hence, \eqref{chiP2} holds.
	
\end{proof}
\begin{proof}[Proof of Theorem~\ref{thm: quantrandcol}]
	For any $P\in\cl{P}$, 
	\begin{align}
	\label{eq: SumOf Two}
	\bb{P}\left[\rchi(P)>\Delta_P\right]&=\bb{P}\left[\tr^2(\rchi P)+\tr\big(P-(\rchi P)^2\big)>\Delta_P^2\right]\nonumber\\
	&\leq\bb{P}\left[\big|\tr(\rchi P)\big|>\frac{\Delta_P}{\sqrt{2}}\right]+\bb{P}\left[\tr\big(P-(\rchi P)^2\big)>\frac{\Delta_P^2}{2}\right].
	\end{align}
	We show that for each $P\in\cl{P}$, both $\bb{P}\Big[|\tr(\rchi P)|>\frac{\Delta_P}{\sqrt{2}}\Big]$ and $\bb{P}\Big[\tr\big(P-(\rchi P)^2\big)>\frac{\Delta_P^2}{2}\Big]$ are at most $1/(4M)$. Then, because of \eqref{eq: SumOf Two} we obtain 
	\begin{align*}
	\bb{P}\left[\bigcap_{P\in\mathcal P}\left\{\rchi(P)\leq\Delta_P\right\}\right]&=1-\bb{P}\left[\bigcup_{P\in\mathcal P}\left\{\rchi(P)>\Delta_P\right\}\right]\\
	&\geq 1-\sum_{P\in\mathcal P}\bb{P}\big[\rchi(P)>\Delta_P\big]\tag{the union bound}\\
	&\geq 1-M\frac{1}{2M}=\frac{1}{2}\cc
	\end{align*}
	which is the first claim of the theorem.\\
	If we put $p_1(x_1,x_2,x_3,x_4)= x_1x_3x_2x_4$ and $p_2(x_1,x_2,x_3,x_4)= x_4-(x_1x_3x_2x_4)^2$, then for a random balanced quantum coloring $\rchi=UDU^*$,
	\begin{align*}
	&f_1(U)\coloneqq\tr(p_1\big(U,U^*,D,P)\big)=\tr(\rchi P),\\
	&f_2(U)\coloneqq\tr(p_1\big(U,U^*,D,P)\big)=\tr\big(P-(\rchi P)^2\big).
	\end{align*}
	Since the singular values of $D$, and each $P$ are at most $1$, by Proposition~\ref{prop: ConcProp} we have positive constants $c_1=c_1(p_1)$, $c_2=c_2(p_2)$ and a natural number $N_0$ such that for any $N>N_0$, $\delta>0$, and $P\in\cl P$
	\begin{align}
	&\bb{P}\Bigg[\Big|\tr(\rchi P)-\bb{E}\left[\tr(\rchi P)\right]\Big|>\delta N\Bigg]\leq2\exp\left(-c_1N^2\delta^2\right)\label{rel: FirstConc}\cc\\
	&\bb{P}\Bigg[\Big|\tr\big(P-(\rchi P)^2\big)-\bb{E}\left[\tr(P-(\rchi P)^2)\right]\Big|>\delta N\Bigg]\leq2\exp\left(-c_2N^2\delta^2\right)\label{rel: SecConc}\cdot
	\end{align}
	
	Let $c = \min(c_1,c_2)$. Since the function $x\mapsto\exp(-x)$ is decreasing, $c_1$, $c_2$, and so $c$ are assumed, without loss of generality, to be in the interval $(0,1]$. 
	
	\vspace{.4cm}\noindent
	$\boldsymbol{\bb{P}\Big[|\tr(\rchi P)|>\frac{\Delta_P}{\sqrt{2}}\Big]\leq\frac{1}{4M}:}$\vspace{.1cm}
	
	\noindent
	Note that $\frac{N^2}{N^2-1}\ra(P)- \frac{N}{N^2-1}\ra^2(P)\geq 0$.
	\begin{align*}
	\bb{P}\left[|\tr(\rchi P)|>\frac{\Delta_P}{\sqrt{2}}\right]&\leq\bb{P}\left[|\tr(\rchi P)|>\sqrt{\frac{1}{c}\log(8M)}+\frac{\ra(P)}{N}\right]\\
	&\leq\bb{P}\left[\big|\tr(\rchi P)\big|-\frac{\ra(P)}{N}>\sqrt{\frac{1}{c_1}\log(8M)}\right]\tag{$c\leq c_1$}\\
	&\leq\bb{P}\left[\big|\tr(\rchi P)\big|-\big|\bb E\big[\tr(\rchi P)\big]\big|>\sqrt{\frac{1}{c_1}\log(8M)}\right]\tag{by \eqref{chiP}}\\
	&\leq\bb{P}\left[\Big|\tr(\rchi P)-\bb E\big[\tr(\rchi P)\big]\Big|>\sqrt{\frac{1}{c_1}\log(8M)}\right]\\
	&\leq2\exp\left[-\log(8M)\right]\tag{for large values of $N$, by\eqref{rel: FirstConc}}\\
	&=\frac{1}{4M}\cdot
	\end{align*}
	$\boldsymbol{\bb{P}\left[\tr\big(P-(\rchi P)^2\big)>\frac{\Delta_P^2}{2}\right]\leq\frac{1}{4M}:}$
	\begin{align*}
	\bb{P}&\left[\tr\big(P-(\rchi P)^2\big)>\frac{\Delta_P^2}{2}\right]\\
	&\leq\bb{P}\left[\tr\big(P-(\rchi P)^2\big)>\frac{1}{c}\log(8M)+\frac{N^2}{N^2-1}\ra(P)-\frac{N}{N^2-1}\ra^2(P)\right]\\
	&\leq\bb{P}\left[\tr\big(P-(\rchi P)^2\big)-\bb{E}\big[\tr\big(P-(\rchi P)^2\big)\big]>\frac{1}{c_2}\log(8M)\right]\tag{{$c\leq c_2$, \eqref{chiP2}}}\\
	&\leq\bb{P}\left[\Big|\tr\big(P-(\rchi P)^2\big)-\bb{E}\big[\tr\big(P-(\rchi P)^2\big)\big]\Big|>\sqrt{\frac{1}{c_2}\log(8M)}\right]\tag{$c_2\leq1$}\\
	&\leq\frac{1}{4M}\cdot\tag{by \eqref{rel: SecConc}}
	\end{align*}
To prove the second claim, note that for large values of $N$, $\bb{P}\left[\bigcap_{P\in\mathcal P}\left\{\rchi(P)\leq\Delta_P\right\}\right]>0$, so there exists a realization $\wt\rchi$ of $\rchi$ for which 
		\begin{align*}
		\max_{P\in\cl{P}}\wt{\rchi}(P)\leq\max_{P\in\cl{P}}\Delta_P,
		\end{align*}
and thus, $\qdisc(\cl{P})\leq\max_{P\in\cl{P}}\Delta_P$. Moreover, since $0\leq\ra(P)\leq N$, for every $P\in\cl{P}$,
	\begin{align*}
	\Delta_P&\leq\sqrt{2}\left[\sqrt{\frac{1}{c}\log(8M)+\frac{N^3}{N^2-1}}+1\right]\\
	&\leq\sqrt{2}\left[\sqrt{\frac{1}{c}\log(8M)+N+1}+1\right].\tag{$N\geq 2$}
	\end{align*}
Hence, $\qdisc(P)=O(\sqrt{N+\log M})$.
	
\end{proof}
\subsection{Proof of Theorem~\ref{Thm: Tightness}}
A probabilistic proof of Theorem~\ref{Thm: Tightness} is given in this subsection. It is based on the behavior of a random projection system. Note that for $\Pi\coloneqq\left[\begin{array}{c|c}
I_{\lfloor\frac{N}{2}\rfloor}&\mathbf{{0}}\\
\hline
\mathbf{0}&\mathbf{0}\end{array}\right]\in\cl M_N(\bb C)$, and a Haar-distributed random element $U$ of $\bb{U}(N)$, the random matrix $P\coloneqq U\Pi U^*$ is a random orthogonal projection. In fact, $P$ is the orthogonal projection onto a random subspace of dimension $\lfloor N/2\rfloor$ in $\bb{C}^N$. In this subsection we assume that $\cl P = \left\{P_1, \ldots, P_M\right\}$, where for each $i$, $P_i = U_i\Pi U_i^*$ and $\{U_1, \ldots, U_M\}$ is an independent set of Haar-distributed random elements of $\bb{U}(N)$. Moreover, we suppose that $M$ and $N$ satisfy the conditions of Theorem~\ref{Thm: Tightness}, i.e. for some $a,b>0$
\begin{align}
\label{rel: AcceptRegimes}
aN\leq M\quad,\quad \log M\leq bN.
\end{align} 
To prove the theorem, it suffices to show that there exists a constant $c=c(b)>0$ for which
\begin{align}
\label{rel: ProbLessOne}
	\bb{P}\left[ \qdisc(\cl{P})\leq c\sqrt{N+\log M}\right]<1.
\end{align}
Note that for the $\rchi(P)$ given by \eqref{eq: ChiP} and any $c>0$, 
\begin{align*}
	\bb{P}\left[ \qdisc(\cl{P})\leq c\sqrt{N+\log M}\right]&=\bb{P}\left[ \min_{\rchi}\max_{P\in\cl{P}}\rchi(P)\leq c\sqrt{N+\log M}\right]\\
	&\leq\bb{P}\left[ \bigcup_{\rchi}\left(\max_{P\in\cl{P}}\rchi(P)\leq c\sqrt{N+\log M}\right)\right]\cdot
\end{align*}
If the number of quantum colorings were finite and not too large, using the union bound together with Proposition~\ref{prop: ConcProp} could be helpful to obtain \eqref{rel: ProbLessOne}. Although the set of quantum colorings is uncountable, the same idea can be applied in an indirect way. The method is borrowed from Subsection 2.3.1 of~\cite{TaoMatrixBook}. There, an upper bound for the operator norm of a certain family of random matrices is provided using this technique, which is named the epsilon-net argument. 
The following lemma provides a way to go from the uncountable set of quantum colorings to an appropriate finite subset.
\begin{lemma}
	\label{lem: EpsilonNet}
	Suppose $C\subseteq\bb{R}^d$ is compact. For any $\epsilon>0$, there exists a finite set $\Sigma^{^{(\epsilon)}}\subseteq C$ so that for each $L$-Lipschitz function $\varphi:C\to\bb{R}$, and $\Delta\in\bb{R}$, if $\min_{x\in C}\varphi(x)\leq\Delta$, then
	$\min_{x\in\Sigma^{^{(\epsilon)}}}\varphi(x)\leq\Delta+\epsilon L$.
\end{lemma} 
\begin{proof}
	Set $\Sigma^{^{(\epsilon)}}$ to be an $\epsilon$-net in $C$, that is a maximal subset of $C$ with the property that the distance between each two distinct points is more than $\epsilon$. \\
	To see that $\Sigma^{^{(\epsilon)}}$ is finite,	suppose $B(x,r)$ and $\bar{B}(x,r)$ are open and closed balls with radius $r>0$ and center $x\in\bb{R}^d$, respectively. By compactness, $C$ is included in $B(O,r)$ for some $r>0$ ($O$ is the origin). For any two distinct $x,y\in\Sigma^{^{(\epsilon)}}$, $\bar{B}(x,\epsilon/2)$ and $\bar{B}(y,\epsilon/2)$ are disjoint. Also, each of them is included in $\bar{B}(O,r+\epsilon/2)$. Thus,
	\begin{align*}
		&\sum_{x\in\Sigma^{^{(\epsilon)}}}\vol(\bar{B}(x,\frac\epsilon 2))=\vol(\bigcup_{x\in\Sigma^{^{(\epsilon)}}}\bar{B}(x,\frac\epsilon 2))\leq\vol(\bar{B}(O,r+\frac\epsilon 2)).
	\end{align*}
		Hence,
	\begin{align*}
		&\left|\Sigma^{^{(\epsilon)}}\right|\leq \frac{\vol(\bar{B}(O,r+\frac\epsilon 2))}{\vol(\bar{B}(O,\frac\epsilon 2))}=\left(\frac{r+\frac{\epsilon}{2}}{\frac\epsilon 2}\right)^d=\left(\frac{2r+\epsilon}{\epsilon}\right)^d<\infty.
	\end{align*}
	For a function $\varphi$ satisfying the given conditions, if $\min_{x\in C}\varphi(x)\leq\Delta$, there is some $x_0\in C$ such that $\varphi(x_0)\leq\Delta$. Because of maximality of $\Sigma^{^{(\epsilon)}}$, $\|x_1-x_0\|\leq\epsilon$ for at least one point $x_1\in\Sigma^{^{(\epsilon)}}$. By the Lipschitz property,
	\begin{align*}
		&\varphi(x_1)-\varphi(x_0)\leq|\varphi(x_1)-\varphi(x_0)|\leq L\norm{x_1-x_0}{}.
	\end{align*}
	Thus,  $\varphi(x_1)\leq\Delta+\epsilon L$, and hence $\min_{x\in\Sigma^{^{(\epsilon)}}}\varphi(x)\leq\Delta+\epsilon L$.

\end{proof}
\begin{corollary}
	\label{Cor: RandLipFunc}
	Suppose $Y$ is a random object and the random function $\varphi(x,Y)$ has the property that for any realization $y$ of $Y$ its value, $\varphi(x;y)$, is an $L$-Lipschitz real-valued function on a compact set $C$. Then, for any $\Delta,\epsilon>0$ 
	\begin{align}
	\label{rel: RandLipFunc}
		\bb{P}\left[\min_{x\in C}\varphi(x,Y)\leq\Delta\right]\leq\sum_{x\in\Sigma^{^{(\epsilon)}}}\bb{P}\big[\varphi(x,Y)\leq\Delta+\epsilon L\big]\cdot
	\end{align} 
\end{corollary}
\begin{proof}
 By Lemma~\ref{lem: EpsilonNet}, the same $\Sigma^{^{(\epsilon)}}$ works for each  $L$-Lipschitz function $\varphi(x;y)$. Moreover, for any realization $y$ of $Y$, if $\min_{x\in C}\varphi(x;y)\leq\Delta$, then $\min_{x\in \Sigma^{^{(\epsilon)}}}\varphi(x;y)\leq\Delta+\epsilon L$. Hence,
 
 \begin{align*}
 \bb{P}\left[\min_{x\in C}\varphi(x,Y)\leq\Delta\right]&\leq\bb{P}\left[\min_{x\in \Sigma^{^{(\epsilon)}}}\varphi(x,Y)\leq\Delta+\epsilon L\right]\\
 &=\bb{P}\left[\bigcup_{x\in \Sigma^{^{(\epsilon)}}}\left\{\varphi(x,Y)\leq\Delta+\epsilon L\right\}\right]\\
 &\leq\sum_{x\in\Sigma^{^{(\epsilon)}}}\bb{P}\big[\varphi(x,Y)\leq\Delta+\epsilon L\big]\cdot\tag{the union bound}
 \end{align*} 
 
\end{proof}

Note that appropriate upper bounds on the size of $\Sigma^{^{(\epsilon)}}$, and values of the summands on the right-hand side of \eqref{rel: RandLipFunc}, would result in controlling the probability of the event $\min_{x\in C}\varphi(x,Y)\leq\Delta$.

To verify the Lipschitz property of the functions within the proof, a matrix version of H\"{o}lder inequality will be used.
\begin{proposition}[Tracial H\"{o}lder Inequality]
	For $P,Q\in\cl{M}_N(\bb{C})$ and $1\leq p,q\leq\infty$ with the property that $1/p+1/q=1$, 
	\begin{align}
		\label{rel: Holdeq}
		|\tr(P^*Q)|\leq\norm{P}{p}\norm{Q}{q}\cdot
	\end{align}
	Here, $\norm{.}{p}$ is the Schatten $p$-norm which is defined for $1\leq p\leq\infty$ and $P\in \cl{M}_N(\bb{C})$ with singular values $\sigma_1\leq\sigma_2\leq\cdots\leq\sigma_N$ as
	\begin{align*}
		\norm{P}{p}&=\big[\sum_i\sigma_i^p\big]^{1/p}\quad \text{for}\quad 1\leq p<\infty\quad,\quad\norm{P}{\infty}= \sigma_N\cdot
	\end{align*}
\end{proposition}
This inequality proved in~\cite{Carlen5} is due to Carlen.
\begin{proof}[Proof of Theorem~\ref{Thm: Tightness}]
	Suppose $c>0$ and set $C\subseteq\cl{M}_N(\bb C)$ to be the set of quantum colorings. 
	Put
	\begin{align*}
		C_1=\{\rchi\in{\raisebox{-.87pt}{\scalebox{1}{$C$}}}: |\tr(\rchi)|>\sqrt{N}\log N\}\quad&,\quad C_2=\{\rchi\in {\raisebox{-.87pt}{\scalebox{1}{$C$}}}: |\tr(\rchi)|\leq\sqrt{N}\log N\}\cdot
	\end{align*} 
	Then,
	\begin{align*}
		\bb{P}&\left[\qdisc(\cl P)\leq c\sqrt{N+\log M}\right]\\
		&=\bb{P}\left[\min_{\lrchi}\max_{P\in\cl{P}}\left[\tr^2(\rchi P)+\tr\big(P-(\rchi P)^2\big)\right]\leq c^2(N+\log M)\right]\\
		&\leq\bb{P}\left[\min_{\lrchi\in {\raisebox{-.87pt}{\scalebox{.7}{$C_1$}}}}\max_{P\in\cl{P}}\left[\tr^2(\rchi P)+\tr\big(P-(\rchi P)^2\big)\right]\leq c^2(N+\log M)\right]\\
		&\quad+\bb{P}\left[\min_{\lrchi\in {\raisebox{-.87pt}{\scalebox{.7}{$C_2$}}}}\max_{P\in\cl{P}}\left[\tr^2(\rchi P)+\tr\big(P-(\rchi P)^2\big)\right]\leq c^2(N+\log M)\right]\cdot
	\end{align*}
	 \noindent
	Because of the non-negativity of $\tr^2(\rchi P)$ and $\tr\left(P-(\rchi P)^2\right)$, 
	\begin{align*}
		\bb{P}\left[\qdisc(\cl P)\leq c\sqrt{N+\log M}\right]&\leq\bb{P}\left[\min_{\lrchi\in {\raisebox{-.87pt}{\scalebox{.7}{$C_1$}}}}\max_{P\in\cl{P}}\big|\tr(\rchi P)\big|\leq c\sqrt{N+\log M}\right]\\
		&\quad+\bb{P}\left[\min_{\lrchi\in {\raisebox{-.87pt}{\scalebox{.7}{$C_2$}}}}\max_{P\in\cl{P}}\tr\big(P-(\rchi P)^2\big)\leq c^2(N+\log M)\right]\cdot
	\end{align*}
	It is enough to prove the existence of a constant $ c>0$ which depends only on $b$, given in \eqref{rel: AcceptRegimes}, and satisfies the following inequalities for all but a finite set of values of $(M,N)$. 
	\begin{align*}
		&\bb{P}\left[\min_{\lrchi\in {\raisebox{-.87pt}{\scalebox{.7}{$C_1$}}}}\max_{P\in\cl{P}}\big|\tr(\rchi P)\big|\leq c\sqrt{N+\log M}\right]<\frac{1}{2}\cc\\\ &\bb{P}\left[\min_{\lrchi\in {\raisebox{-.87pt}{\scalebox{.7}{$C_2$}}}}\max_{P\in\cl{P}}\tr\big(P-(\rchi P)^2\big)\leq c^2(N+\log M)\right]<\frac{1}{2}\cdot
	\end{align*}
To achieve this goal  the epsilon-net argument will be applied twice. To be permitted to make use of Lemma~\ref{lem: EpsilonNet}, the compactness and Lipschitz properties must be checked.\\
The space $\left(\cl{M}_N(\bb{C}),\norm{.}{2}\right)$ is isomorphic to $(\bb{R}^{2N^2},\norm{.}{2})$ as a normed vector space. The subset of Hermitian matrices, $\cl{H}_N$, is actually an $N^2$-dimensional subspace of $\cl{M}_N(\bb{C})$, and so a Euclidean metric space. Eigenvalues of a matrix are continuous functions of its entries. Moreover, $\norm{\rchi}{2}= \sqrt{N}$ and  $\tr(\rchi)\in\{\pm 1,\pm 2, \ldots, \pm N\}$ for each $\rchi\in C$. Hence, $C_1$ and $C_2$ are compact sets in $\cl{H}_N$.\\
The maximum of a family of $L$-Lipschitz functions is again  $L$-Lipschitz. Hence, according to Corollary~\ref{Cor: RandLipFunc}, it suffices to verify the Lipschitz property of the functions $\big|\tr(\rchi \tilde{P})\big|$ and
			$\tr\left(\tilde{P}-(\rchi \tilde P)^2\right)$ for $\tilde{P}$, a realization of  $P\in\cl P$.
			\begin{align*}
				\Big|\big|\tr(\rchi_1 \tilde{P})\big|-\big|\tr(\rchi_2 \tilde{P})\big|\Big|
				&\leq\big|\tr(\rchi_1-\rchi_2)\tilde{P}\big|\\
				&\leq\norm{\tilde{P}}{2}\norm{\rchi_1-\rchi_2}{2}\tag{by \eqref{rel: Holdeq}}\\
				&\leq\sqrt{\frac N 2}\norm{\rchi_1-\rchi_2}{2}.\tag{$\ra(\tilde{P}) =\lfloor\frac N 2\rfloor$}
			\end{align*}
			Hence, for any realization $\wt{\cl{P}}$ of $\cl{P}$, $\varphi_1(\rchi;\wt{\cl{P}})\coloneqq\max_{\tilde{P}\in\tilde{\cl{P}}}\big|\tr(\rchi \tilde{P})\big|$ is $\sqrt{N/2}$-Lipschitz on $C_1$. 
			On the other hand, since $\|{\tilde{P}}\|_{\infty}=1$, and $\norm{\rchi_1}{1}=\norm{\rchi_2}{1}=N$,
			
			\begin{align*}
				\left|\tr\!\left(\!\tilde{P}-(\rchi_1 \tilde P)^2\!\right)\!-\!\tr\!\left(\!\tilde{P}-(\rchi_2 \tilde P)^2\!\right)\!\right|
				&=\big|\tr(\rchi_2\tilde{P}\rchi_2\tilde{P}-\rchi_2\tilde{P}\rchi_1 \tilde{P}+\rchi_2\tilde{P}\rchi_1 \tilde{P}-\rchi_1 \tilde{P} \rchi_1 \tilde{P})\big|\\
				&=\big|\tr[(\rchi_2+\rchi_1)\tilde{P}(\rchi_2-\rchi_1)\tilde{P}]\big|\\
				&\leq\|{\tilde{P}(\rchi_2-\rchi_1)\tilde{P}}\|_{\infty}\norm{\rchi_2+\rchi_1}{1}\tag{by \eqref{rel: Holdeq}}\\
				&\leq\|{\tilde{P}}\|_{\infty}^2\norm{\rchi_2-\rchi_1}{\infty}\big(\norm{\rchi_1}{1}+\norm{\rchi_2}{1}\big)\\
				&\leq 2N\norm{\rchi_2-\rchi_1}{\infty}\\
				&\leq 2N\norm{\rchi_2-\rchi_1}{2}\cdot
			\end{align*}
			Therefore, for each realization $\wt{\cl{P}}$, the function $\varphi_2(\rchi;\wt{\cl{P}})\coloneqq\max_{\tilde{P}\in\wt{\cl{P}}}\tr\left(\tilde{P}-(\rchi \tilde P)^2\right)$ is $2N$-Lipschitz on $C_2$. \\
	We conclude that for any $\epsilon>0$ there are (finite) $\epsilon$-nets $\Sigma_1^{^{(\epsilon)}}\subseteq C_1$ and $\Sigma_2^{^{(\epsilon)}}\subseteq C_2$ such that
	\begin{align}
		\label{rel: InfToFin1}
		&\bb{P}\left[\min_{\lrchi\in {\raisebox{-.87pt}{\scalebox{.7}{$C_1$}}}}\max_{P\in\cl{P}}\big|\tr(\rchi P)\big|\leq c\sqrt{N+\log M}\right]\nonumber\\
		&\quad\leq\bb{P}\left[\min_{\lrchi\in {\raisebox{-.87pt}{\scalebox{.7}{$\Sigma_1^{^{(\epsilon)}}$}}}}\max_{P\in\cl{P}}\left|\tr(\rchi P)\right|\leq c\sqrt{N+\log M}+\epsilon\sqrt{\frac{N}{2}}\right]\nonumber\\
		&\quad\leq\sum_{\lrchi\in {\raisebox{-.87pt}{\scalebox{.7}{$\Sigma_1^{^{(\epsilon)}}$}}}}\!\bb{P}\!\left[\max_{P\in\cl{P}}\left|\tr(\rchi P)\right|\leq c\sqrt{N+\log M}+\epsilon\sqrt{\frac{N}{2}}\right],
	\end{align}
	and
	\begin{align}
		\label{rel: InfToFin2}
		&\bb{P}\left[\min_{\lrchi\in {\raisebox{-.87pt}{\scalebox{.7}{$C_2$}}}}\max_{P\in\cl{P}}\tr\big(P-(\rchi P)^2\big)\leq c^2(N+\log M)\right]\nonumber\\
		&\quad\leq\sum_{\lrchi\in{\raisebox{-.87pt}{\scalebox{.7}{$\Sigma_2^{^{(\epsilon)}}$}}}}\bb{P}\left[\max_{P\in\cl{P}}\tr\left(P-(\rchi P)^2\right)\leq c^2(N+\log M)+2\epsilon N \right]\cdot
	\end{align}
	We make two observations. First suppose $0<\epsilon\leq 1$. According to the proof of Lemma~\ref{lem: EpsilonNet}, and since for every quantum coloring $\rchi$, $\norm{\rchi}{2}=\sqrt{N}$, 
			\begin{align}
				\label{rel: SizeOfNets}
				\left|\Sigma_1^{^{(\epsilon)}}\right|,\left|\Sigma_2^{^{(\epsilon)}}\right|\leq\left( \frac{2\sqrt{N}+\epsilon}{\epsilon}\right)^{N^2}\leq\left( \frac{3\sqrt{N}}{\epsilon}\right)^{N^2}.
			\end{align}
	 On the other hand, if $\Delta >0$, then for any quantum coloring $\rchi$
			\begin{align}
				\bb{P}\left[\max_{P\in\cl{P}}\left|\tr(\rchi P)\right|\leq\Delta\right]& = \bb{P}\left[\bigcap_{P\in\cl{P}}\left\{\left|\tr(\rchi P)\right|\leq\Delta\right\}\right]\nonumber \\
				& = \prod_{P\in\cl{P}}\bb{P}\left[\left|\tr(\rchi P)\right|\leq\Delta\right]\tag{by independence}\nonumber\\
				\label{rel: MaxToP1}
				& = \left(\bb{P}\left[\left|\tr(\rchi P)\right|\leq\Delta\right]\right)^M,
			\end{align}
			where $P$ is an arbitrary element of $\cl P$. The last equality holds because elements of $\cl{P}$ are identically distributed.\\
			With the same argument, for any quantum coloring $\rchi$ 
			\begin{align}
				\label{rel: MaxToP2}
				\bb{P}\left[\max_{P\in\cl{P}}\tr\left(P-(\rchi P)^2\right)\leq\Delta\right]= \left(\bb{P}\left[\tr\left(P-(\rchi P)^2\right)\leq\Delta\right]\right)^M.
		\end{align}
	Combining \eqref{rel: InfToFin1}, \eqref{rel: SizeOfNets}, and \eqref{rel: MaxToP1} we get
	\begin{align}
		\label{rel: BeforeConc1}
		&\bb{P}\left[\min_{\lrchi\in{\raisebox{-.87pt}{\scalebox{.7}{$C_1$}}}}\max_{P\in\cl{P}}\big|\tr(\rchi P)\big|\leq c\sqrt{N+\log M}\right]\nonumber\\
		&\quad\leq\left(\frac{3\sqrt{N}}{\epsilon}\right)^{N^2}\left(\max_{\lrchi\in  {\raisebox{-.87pt}{\scalebox{.7}{$\Sigma_1^{^{(\epsilon)}}$}}}}\bb{P}\left[\left|\tr(\rchi P)\right|\leq c\sqrt{N+\log M}+\epsilon\sqrt{\frac{N}{2}}\right]\right)^M.
	\end{align}
	Similarly, \eqref{rel: InfToFin2}, \eqref{rel: SizeOfNets}, and \eqref{rel: MaxToP2} lead to
	\begin{align}
		\label{rel: BeforeConc2}
		&\bb{P}\left[\min_{\lrchi\in{\raisebox{-.87pt}{\scalebox{.7}{$C_2$}}}}\max_{P\in\cl{P}}\tr\big(P-(\rchi P)^2\big)\leq c^2(N+\log M)\right]\nonumber\\
		&\quad\leq\left(\frac{3\sqrt{N}}{\epsilon}\right)^{N^2}\left(\max_{\lrchi\in{\raisebox{-.87pt}{\scalebox{.7}{$\Sigma_2^{^{(\epsilon)}}$}}}}\bb{P}\left[ \tr\big(P-(\rchi P)^2\big)\leq c^2(N+\log M)+2\epsilon N\right]\right)^M.
	\end{align}
	It remains to provide upper bounds for $\max_{\lrchi\in {\raisebox{-.87pt}{\scalebox{.7}{$\Sigma_1^{^{(\epsilon)}}$}}}}\bb{P}\left[\left|\tr(\rchi P)\right|\leq c\sqrt{N+\log M}+\epsilon\sqrt{N/2}\right]$ and $\max_{\lrchi\in \Sigma_2^{^{(\epsilon)}}}\bb{P}\left[ \tr\big(P-(\rchi P)^2\big)\leq c^2(N+\log M)+2\epsilon N\right]$, and then to find proper values of $\epsilon,  c>0$. 
	Calculation of these bounds will be done by the concentration inequality given in Proposition~\ref{prop: ConcProp}. First, we bound the expected values. Considering $\rchi=VDV^*$ as a spectral decomposition of $\rchi$, it holds that
	\begin{align*}
		\tr(\rchi P)\stackrel{d}{=}\tr(UDU^*\Pi)\quad,\quad\tr\left(P-(\rchi P)^2\right)\stackrel{d}{=}\tr\left(\Pi-(UDU^*\Pi)^2\right)
	\end{align*} 
	for a random unitary $U$ with Haar distribution. Therefore, using the calculation done in the last subsection, we obtain that for $\rchi\in C_1$, 
	\begin{align}
		\label{rel: ExpEst1}
		\left|\bb{E}\left[\tr(\rchi P)\right]\right|=\left|\frac{1}{N}\tr(\rchi)\tr(\Pi)\right|\stackrel{N\geq 3}{\geq}\frac 1 3\sqrt{N}\log N\cc
	\end{align} 
	and for $\rchi\in C_2$,
	\begin{align}
		\label{rel: ExpEst2}
		\bb{E}\left[\tr\left(P-(\rchi P)^2\right)\right]&=-\tr^2(\rchi)\lfloor\frac N 2\rfloor\frac{N-\lfloor\frac N 2\rfloor}{N(N^2-1)}+N\lfloor\frac N 2\rfloor\frac{N-\lfloor\frac N 2\rfloor}{N^2-1}\nonumber\\
		&\geq \lfloor\frac N 2\rfloor\frac{N-\lfloor\frac N 2\rfloor}{N^2-1}\left(-\log^2(N)+N\right)\nonumber\\
		&\stackrel{N\geq 152}{\geq}\frac{N}{6}\cdot
	\end{align}
	Now, we calculate the upper bounds.
	\begin{align}
		\hspace{-.4cm}\bb{P}&\left[\big|\tr(\rchi P)\big|\leq c\sqrt{N+\log M}+\epsilon\sqrt{\frac{N}{2}}\right]\nonumber\\
		\hspace{-.4cm}&\leq\bb{P}\left[\big|\tr(\rchi P)\big|\leq(c+\epsilon)\sqrt{N+\log M}\right]\nonumber\\
    	\hspace{-.4cm}&\leq\bb{P}\left[\big|\tr(\rchi P)\big|-\big|\bb{E}\left[\tr(\rchi P)\right]\big|\leq(c+\epsilon)\sqrt{1+b}\sqrt{N}-\frac{1}{3}\sqrt{N}\log N\right]\tag{\eqref{rel: AcceptRegimes}, \eqref{rel: ExpEst1}, $N\geq 3$}\nonumber\\
		\hspace{-.4cm}&\leq\bb{P}\left[\big|\tr(\rchi P)\big|-\big|\bb{E}\left[\tr(\rchi P)\right]\big|\leq-\frac{1}{4}\sqrt{N}\log N\right]\tag{for large $N$}\nonumber\\
		\hspace{-.4cm}&\leq\bb{P}\left[\big|\tr(\rchi P)-\bb{E}\left[\tr(\rchi P)\right]\big|\geq\frac{1}{4}\sqrt{N}\log N\right]\nonumber\\
		\label{rel: UpperBound1}
		\hspace{-.4cm}&\leq2\exp\left[-\frac{c_1}{16}N\log^2N\right]\cdot
	\end{align}
	\noindent
	The inequality \eqref{rel: UpperBound1} is true by Proposition~\ref{prop: ConcProp} for large enough $N$s and the constant $c_1$ introduced in the previous subsection.
\noindent
\begin{align}
		\bb{P}&\left[\tr\big(P-(\rchi P)^2\big)\leq c^2(N+\log M)+2\epsilon N\right]\nonumber\\
		&\leq\bb{P}\Big[\tr\big(P-(\rchi P)^2\big)-\bb{E}\big[\tr\big(P-(\rchi P)^2\big)\big]\leq c^2(N+\log M)+2\epsilon N-\frac{N}{6}\Big]\tag{\scriptsize{\eqref{rel: ExpEst2}, $N\geq 152$}}\nonumber\\
		&\leq\bb{P}\Big[\tr\big(P-(\rchi P)^2\big)-\bb{E}\big[\tr\big(P-(\rchi P)^2\big)\big]\leq\big(c^2(1+b) +2\epsilon -\frac{1}{6}\big)N\Big]\tag{\eqref{rel: AcceptRegimes}, $N\geq 152$}\nonumber\\
		&\leq\bb{P}\left[\Big|\tr\big(P-(\rchi P)^2\big)-\bb{E}\big[\tr\big(P-(\rchi P)^2\big)\big]\Big|\geq \big(\frac{1}{6}- c^2(1+b)-2\epsilon\big)N\right]\cdot\nonumber
	\end{align}
	The last inequality is valid if 
	\begin{align}
		\label{rel: Cond1}
		\frac{1}{6}- c^2(1+b)-2\epsilon>0\cdot
	\end{align} 
	Assuming \eqref{rel: Cond1}, we conclude from Proposition~\ref{prop: ConcProp} that
	\begin{align}
		\label{rel: UpperBound2}
		&\hspace{-.1cm}\bb{P}\left[\tr\big(P-(\rchi P)^2\big)\leq c^2(N+\log M)+2\epsilon N\right]\leq 2\exp\left[-c_2\big(\frac{1}{6}- c^2(1+b)-2\epsilon\big)^2N^2 \right],
	\end{align}
	where $c_2$ is the same as before and $N$ is large enough.\\
Substituting \eqref{rel: UpperBound1} in \eqref{rel: BeforeConc1} gives
	\begin{align*}
		\bb{P}\left[\min_{\lrchi\in {\raisebox{-.87pt}{\scalebox{.7}{$C_1$}}}}\max_{P\in\cl{P}}\big|\tr(\rchi P)\big|\leq c\sqrt{N+\log M}\right]
		&\leq\left(\frac{3\sqrt{N}}{\epsilon}\right)^{N^2}\left(2\exp\left[-\frac{c_1}{16}N\log^2N\right]\right)^M\\
		&=\exp\left[N^2\log\frac{3}{\epsilon}+\frac{N^2}{2}\log N+M\log 2\right.\\
		&\qquad\qquad\left.-\frac{c_1}{16}NM\log^2N\right]\\
		&\leq\exp\left[N^2\log N-\frac{c_1}{17}NM\log^2N \right]\tag{\small{for large $N$}}\\
		&\leq\exp\left[N^2\log N-\frac{c_1}{17}aN^2\log^2N \right]\tag{by \eqref{rel: AcceptRegimes}}\\
		&\to 0\cdot\tag{as $N\to\infty$}
	\end{align*}
	Hence, $\bb{P}\left[\min_{\lrchi\in {\raisebox{-.87pt}{\scalebox{.7}{$C_1$}}}}\max_{P\in\cl{P}}\big|\tr(\rchi P)\big|\leq c\sqrt{N+\log M}\right]$ is less than $1/2$ for large enough values of $N$.\\
	Similarly, \eqref{rel: UpperBound2} and \eqref{rel: BeforeConc2} give us that provided \eqref{rel: Cond1},
	\begin{align*}
		\bb{P}&\left[\min_{\lrchi\in {\raisebox{-.87pt}{\scalebox{.7}{$C_2$}}}}\max_{P\in\cl{P}}\tr\big(P-(\rchi P)^2\big)\leq c^2(N+\log M)\right]\\
		&\leq\exp\left[N^2\log N-\frac{c_2}{2}\big(\frac{1}{6}-c^2(1+b)-2\epsilon\big)^2\!N^2 M\right]\tag{for large $N$}\\
		&\leq\exp\left[N^2\log N-\frac{c_2}{2}\big(\frac{1}{6}-c^2(1+b)-2\epsilon\big)^2aN^3\right]\tag{by \eqref{rel: AcceptRegimes}}\\
		&\to 0\cdot\tag{as $N\to\infty$}
	\end{align*}
	Hence, $\bb{P}\left[\min_{\lrchi\in {\raisebox{-.87pt}{\scalebox{.7}{$C_2$}}}}\max_{P\in\cl{P}}\tr\big(P-(\rchi P)^2\big)\leq c^2(N+\log M)\right]$ will, ultimately, be less than $1/2$.\\
	Note that up to this point, no restriction has been imposed on $M$, and $N$ has been restricted finitely many times. Thus, to finish the proof it is enough to show that for any $a,b>0$, \eqref{rel: Cond1} is satisfied by appropriate values of $\epsilon, c=c(b)>0$ which are independent of $M$, $N$. This is true since we can put $\epsilon=1/24$ and $c=1/(4\sqrt{1+b})$. 
	
\end{proof}

\begin{remark}
	The range of tightness which is given in Theorem~\ref{Thm: Tightness}, includes projection systems with a size that can vary from linear to exponential with respect to the dimension of the ground space. Comparing with the combinatorial case, the quantum random coloring upper bound is tight in a wider range. The reason seems to be hidden in the strength of the concentration properties of random unitary matrices, e.g. what explained in Proposition~\ref{prop: ConcProp}. 
\end{remark}


\section{Quantum Discrepancy of Set Systems}
\label{Sec: SetSysProofs}
\subsection{Proof of Theorem~\ref{Thm: ZerosDiscQDisc}}
Proof of Theorem~\ref{Thm: ZerosDiscQDisc} is based on the formulation given in Corollary~\ref{Col: NormFormulation}, i.e.
\begin{align*}
\qdisc(\cl P)=\min_{\rchi}\max_{P\in\cl P}\left[\tr^2(\rchi P)+\frac{1}{2}\|[\rchi,P]\|_2^2\right]^{1/2}.
\end{align*}
\begin{proof}[ Proof of Theorem~\ref{Thm: ZerosDiscQDisc}]
	Since $0\leq\qdisc(\cl{P}_{\cl S})\leq\disc(\cl S)$, if $\disc(\cl S) = 0$, then we have $\qdisc(\cl{P}_{\cl S})=0$. Now,  suppose that $\qdisc(\cl{P}_{\cl S})=0$. By the above formula, there exists a quantum coloring $\rchi$ such that for each $P\in\cl{P}_{\cl S}$, $\|[\rchi,P]\|_2=0$. This means that each element of $\cl{P}_{\cl S}$ commutes with $\rchi$. Thus $\rchi$ and $P$s are simultaneously unitarily diagonalizable. Assume that $U$ is a unitary matrix such that
	\begin{align*}
	\wt{\rchi}=U\rchi U^*, \quad \wt{P}=UPU^*\quad\text{for}\quad P\in\cl{P}_{\cl S},
\end{align*}
are diagonal. Hence, $\wt{\cl{P}}\coloneqq\{\wt{P}: P\in\cl{P}_{\cl S}\}$ is associated to a (unique) set system $\wt{\cl S}$. By the cyclic property of the trace function, $\qdisc(\wt{\cl P})=\qdisc(\cl P_{\cl S})=0$, and since $\wt{\rchi}$ is diagonal, $\disc(\wt{\cl S})=0$. Now, it is enough to prove that $\disc(\cl S)=\disc(\wt{\cl S})$. In order to establish this equality, we show that there exists a permutation $\sigma$ on $[N]$ such that $\wt{\cl S}=\{\sigma(S):S\in\cl S\}$. Then $\diag(\wt{\rchi})\circ\sigma$ will be a coloring of $[N]$ that produces zero discrepancy in $\cl S$.\\
We define an equivalence relation on $[N]$ as
\begin{align*}
m\sim n\Leftrightarrow\forall S\in\cl S[m\in S\Leftrightarrow n\in S].
\end{align*} 
 Suppose that $C_1,\ldots,C_k$ are the equivalence classes of this relation. Since
 \begin{align*}
 P_{C_{i}}=\prod_{\substack{S\in\cl S\\S\cap C_i\neq\emptyset}} P_S\prod_{\substack{S\in\cl S\\S\cap C_i=\emptyset}}(I-P_S),
 \end{align*}
 $UP_{C_i}U^*$ is a diagonal projection for each $i = 1,2,\ldots, k$. Let $D_i\subseteq[N]$ be the corresponding set of this projection. The map $P\mapsto UPU^*$ is rank-preserving, so $|D_i|=|C_i|$. For each $i = 1,2,\ldots, k$ we fix a bijection $\sigma_i: C_i\to D_i$. Now, $\sigma: [N]\to[N]$ is uniquely determined by the set of relations $\sigma\big\arrowvert_{C_i}=\sigma_i$. The map $\sigma$ is a permutation since the $C_i$s are equivalence classes. It follows from the definition of $\wt{\cl P}$ that $\wt{S}\in\wt{\cl S}$ exactly when there is a set $S\in\cl S$ such that $\wt{S}$ is the set associated to the projection $UP_SU^*$. Consider a pair of $S\in\cl S$ and $\wt{S}\in\wt{\cl S}$ with this relation. We have
 \begin{align}
 \label{eq: TheLast}
 S=\bigcup_{i:C_i\cap S\neq\emptyset}C_i\quad,\quad P_S=\sum_{i:C_i\cap S\neq\emptyset}P_{C_i}.
 \end{align}
 Since $P_{\wt S}=UP_SU^*$, it holds that 
 \begin{align}
 P_{\wt S}&=\sum_{i:C_i\cap S\neq\emptyset}UP_{C_i}U^*\tag{by \eqref{eq: TheLast}}\\
 &=\sum_{i:C_i\cap S\neq\emptyset}P_{D_i}\nonumber\\
 &=\sum_{i:C_i\cap S\neq\emptyset}P_{\sigma(C_i)}\nonumber\\
 &=P_{_{\left[\bigcup_{_{C_i\cap S\neq\emptyset}}\sigma(C_i)\right]}},\tag{$\sigma(C_i)$s are disjoint}
 \end{align}
 and therefore,
  \begin{align*}
  \wt{S}=\bigcup_{i:C_i\cap S\neq\emptyset}\sigma(C_i)=\sigma(\bigcup_{i:C_i\cap S\neq\emptyset}C_i)=\sigma(S).
  \end{align*}
\end{proof}	
\subsection{Proof of Theorem~\ref{Thm: QuantCombThm2}}
\noindent
Suppose $\cl{S}$ is a set system on $[N]$ and $\cl P_{\cl{S}}$ is its corresponding projection system. We know form Lemma~\ref{lem: ExpToMat} that

\begin{align*}
\qdisc(\cl{P}_{\cl{S}})=\min_{\substack{{\mathfrak{X}:\, \text{det. proj.}}\\{\text{process}}}}\max_{S\in\cl{S}}\left[\bb{E}\left[\left(2{\mathfrak{X}}(S)-|S|\right)^2\right]\right]^\frac{1}{2}.
\end{align*} 
The idea of the proof of Theorem~\ref{Thm: QuantCombThm2} is that if there is a determinantal projection process $\mathfrak{X}$ for which $\max_{S\in\cl{S}}\left[\bb{E}\left[\left(2{\mathfrak{X}}(S)-|S|\right)^2\right]\right]^{1/2}$ is very small, then by means of the concentration of $\mathfrak{X}(S)$ around its mean for every $S\in\cl{S}$, we can find a deterministic coloring $\rchi$ which makes $\max_{S\in\cl{S}}|\rchi(S)|$ and so $\disc(\cl{S})$ small. The following propositions provide a representation of $\mathfrak{X}(S)$ for each $S\subseteq [N]$, and a concentration inequality applicable to such representations, respectively.
\begin{proposition}
	\label{prop: Bernoullis}
	Suppose $\mathfrak{X}$ is a determinantal process on $[N]$ and $K$ is its Hermitian kernel. Then,
	\begin{enumerate}[label=\textbf{\roman{enumi}.}]
		\item{For each $S\subseteq[N]$, $\mathfrak{X}\cap S$ is a determinantal process on $S$ that its kernel is equal to $K$ restricted to the set of rows and columns in $S$, i.e. the matrix $\left(K_{ij}\right)_{i,j\in S}$.}
		\item{If $\{\lambda_1,\ldots,\lambda_N\}$ is the set of eigenvalues of $K$, then $\mathfrak{X}([N])$ has the same distribution as $\sum_{i=1}^NX_i$, where $X_i$s are independent and $X_i\sim Bernoulli(\lambda_i)$, i=1,\ldots,N. In particular, if $\mathfrak X$ is a projection determinantal process, i.e. if $K$ is an orthogonal projection, all realizations of $\mathfrak{X}$ are of the same size which is $\tr(K)$.}
	\end{enumerate}	
\end{proposition}
\begin{proposition}[Bernstein's Inequality]
	\label{Prop: BernIneq}
	Consider the set of independent random variables $\{X_i: i\geq 1\}$ with $\bb{E}\left[X_i\right]=0$, and a number $K$ with $\left|X_i\right|\leq K$ for each $i$. For any
	$t>0$ we have
	\begin{align*}
	\bb{P}\left[\Bigg|\sum_{i\geq 1}X_i\Bigg|\geq t\right]\leq 2\exp\left[-\min\left(\frac{t^2}{4\sum_{i\geq 1}\bb{E}[X_i^2]}\,\cc\,\frac{t}{2K}\right)\right].
	\end{align*}
\end{proposition}
See Remark 4.2.5 and Theorem 4.5.3 from~\cite{GAF} for Proposition~\ref{prop: Bernoullis}, and Lemma 2.7.1 in \cite{Talagrand} for Proposition~\ref{Prop: BernIneq}.

\begin{proof}[Proof of Theorem~\ref{Thm: QuantCombThm2}]
	For writing the proof more conveniently, we put 
	\begin{align*}
	A = \max\left[\sqrt{\log(2M)}\ \cc\ \frac{\log(2M)}{\qdisc(\cl{P}_{\cl{S}})}\right]
	\end{align*}
	\noindent
	$\max_{S\in\cl{S}}\left[\bb{E}\left[\left(2{\mathfrak{X}}(S)-|S|\right)^2\right]\right]^{1/2}$ is a continuous function of the kernel of $\mathfrak{X}$ and set of orthogonal projections is compact in $\bb{R}^{N^2}$. Hence, there exists a determinantal projection process $\mathfrak{Y}$ for which 
	\begin{align*}
	\qdisc(\cl{P}_{\cl{S}})=\max_{S\in\cl{S}}\left[\bb{E}\left[\left(2{\mathfrak{Y}}(S)-|S|\right)^2\right]\right]^\frac{1}{2}.
	\end{align*}
	We prove that
	\begin{align}
	\label{eq: QuantCombThm2}
	\bb{P}\left[\max_{S\in\cl{S}}\left|2\mathfrak{Y}(S)-|S|\right|\geq\left(2cA+1\right)\qdisc(\cl{P}_{\cl{S}})\right]<1.
	\end{align}
	This results in the existence of some realization $\wt{\mathfrak{Y}}$ for which
	\begin{align*}
	\max_{S\in\cl{S}}\left|2\wt{\mathfrak{Y}}(S)-|S|\right|<\left(2cA+1\right)\qdisc(\cl{P}_{\cl{S}}).
	\end{align*}
	This relation together with $\max_{S\in\cl{S}}\left|2\wt{\mathfrak{Y}}(S)-|S|\right|\geq\disc(\cl{S})$ completes the proof.\\
	In order to prove \eqref{eq: QuantCombThm2}, it's enough, by the union bound, to show that for each $S\in\cl{S}$
	\begin{align*}
	\label{eq: QuantCombThm-22}
	\bb{P}\left[\left|2\mathfrak{Y}(S)-|S|\right|\geq\left(2cA+1\right)\qdisc(\cl{P}_{\cl{S}})\right]<\frac{1}{M}\cdot
	\end{align*}
	Now, for each $S\in\cl{S}$
	\begin{align*}
	&\bb{P}\left[\left|2\mathfrak{Y}(S)-|S|\right|\geq\left(2cA+1\right)\qdisc(\cl{P}_{\cl{S}})\right]\\
	&\hspace{.5cm}=\bb{P}\left[\left|2\mathfrak{Y}(S)-|S|\right|-\max_{S'\in\cl{S}}\left[\bb{E}\left[\left(2\mathfrak{Y}(S')-|S'|\right)^2\right]\right]^\frac{1}{2}\geq 2cA.\qdisc(\cl{P}_{\cl{S}})\right]\\
	&\hspace{.5cm}\leq\bb{P}\left[\left|2\mathfrak{Y}(S)-|S|\right|-\left[\bb{E}\left[\left(2\mathfrak{Y}(S)-|S|\right)^2\right]\right]^\frac{1}{2}\geq 2cA.\qdisc(\cl{P}_{\cl{S}})\right]\\
	&\hspace{.5cm}\leq\bb{P}\left[\left|2\mathfrak{Y}(S)-|S|\right|-\left|\bb{E}\left[2\mathfrak{Y}(S)-|S|\right]\right|\geq 2cA.\qdisc(\cl{P}_{\cl{S}})\right]\\
	&\hspace{.5cm}\leq\bb{P}\left[\left|2\mathfrak{Y}(S)-|S|-\bb{E}\left[2\mathfrak{Y}(S)-|S|\right]\right|\geq 2cA.\qdisc(\cl{P}_{\cl{S}})\right]\\
	&\hspace{.5cm}\leq\bb{P}\Big[\big|\mathfrak{Y}(S)-\bb{E}\left[\mathfrak{Y}(S)\right]\big|\geq cA.\qdisc(\cl{P}_{\cl{S}})\Big]\cdot
	\end{align*}
	By Proposition \ref{prop: Bernoullis}, $\mathfrak{Y}(S)\stackrel{d}{=}\sum_{i=1}^{|S|}Y_i^{(S)}$ where $Y_i^{(S)}$s are independent Bernoulli random variables. Hence, conditions of Proposition \ref{Prop: BernIneq} are satisfied for $\left\{Y_i^{(S)}-\bb{E}\left[Y_i^{(S)}\right]: 1\leq i\leq |S|\right\}$ with $K=1$, and hence
	\begin{align*}
	&\bb{P}\left[\left|2\mathfrak{Y}(S)-|S|\right|\geq\left(2cA
	+1\right)\qdisc(\cl{P}_{\cl{S}})\right]\\
	&\hspace{.5cm}\leq 2\exp\left[-\min\left(\frac{c^2A^2\max_{S'}\bb{E}\left[\left(2\mathfrak{Y}(S')-|S'|\right)^2\right]}{4\sum_{i}\bb{E}\left[\left(Y_i^{(S)}-\bb{E}\left[Y_i^{(S)}\right]\right)^2\right]}\,\cc\,\frac{c}{2}A.\qdisc(\cl{P}_{\cl{S}})\right)\right].
	\end{align*}
	If $y\geq z$, $\min(x,y)\geq\min(x,z)$, for any $x,y,z\in\bb R$. By this and the independence of the $Y_i^{(S)}$s, we obtain
	\begin{align*}
	&\bb{P}\left[\left|2\mathfrak{Y}(S)-|S|\right|\geq\left(2cA+1\right)\qdisc(\cl{P}_{\cl{S}})\right]\\
	&\hspace{.5cm}\leq 2\exp\left[-\min\left(\frac{c^2A^2\bb{E}\left[\left(2\mathfrak{Y}(S)-|S|\right)^2\right]}{\bb{E}\left[\left(2\mathfrak{Y}(S)-\bb{E}\left[2\mathfrak{Y}(S)\right]\right)^2\right]}\,\cc\,\frac{c}{2}A.\qdisc(\cl{P}_{\cl{S}})\right)\right].
	\end{align*} 
	\noindent
	For a random variable $X$, $\bb{E}\left[\left(X-a\right)^2\right]$ takes its minimum value at $a=\bb E\left[X\right]$. Hence,
	\begin{align*}	
	&\bb{P}\left[\left|2\mathfrak{Y}(S)-|S|\right|\geq\left(2cA+1\right)\qdisc(\cl{P}_{\cl{S}})\right]\\
	&\hspace{.5cm}\leq 2\exp\left[-\min\left(c^2A^2\,\cc\,\frac{c}{2}A.\qdisc(\cl{P}_{\cl{S}}\right)\right]\\
	&\hspace{.5cm}< 2\exp\left[-\log(2M)\right]\tag{for each $c>2$}\\
	&\hspace{.5cm}=\frac{1}{M}\cdot
	\end{align*}
	
	\noindent
	To justify the last inequality, note that depending on whether $0<\!\qdisc(\cl{P}_{\cl{S}})<\!\sqrt{\log(2M)}$ or $\qdisc(\cl{P}_{\cl{S}})\geq\sqrt{\log(2M)}$, we have $A = \frac{\log(2M)}{\qdisc(\cl{P}_{\cl{S}})}$ or $A = \sqrt{\log(2M)}$. Therefore, 
	\begin{align*}
    A^2, A.\qdisc(\cl{P}_{\cl{S}})\geq \log(2M)
	\end{align*}
    for every non-zero value of $\qdisc(\cl P_{\cl S})$.
    
 \end{proof}

\section{Concluding Remarks}
\label{Sec: ConcRems}
Quantum discrepancy, as we defined in this paper, is not the first non-commutative version of discrepancy. In \cite{Weaver} Weaver points briefly to the ``interest in noncommutative discrepancy'', and implicitly gives a formulation for the discrepancy of a set of Hermitian matrices as follows.
\noindent
For a set system $\cl{S}$ on $[N]$ with elements indexed as $S_1,S_2,\ldots,S_M$, it can be proved that 
$	\disc(\cl{S})=\min_{\epsilon_1,\ldots,\epsilon_n\in\{\pm 1\}}\norm{\sum_{i=1}^N\epsilon_iv_i}{\infty}$, where for each $i\in[N]$,  $v_i\in\{0,1\}^M$, and $v_i(j)=1$ exactly when $i\in S_j$. If the $v_i$s are permitted to be arbitrary elements of $\bb{C}^M$ instead of being restricted to $0-1$ vectors, then it makes sense to talk about the discrepancy of a set of vectors. Weaver suggests to generalize this new notion to $M\times M$ Hermitian matrices $A_1, \ldots,A_N$ as
\begin{align*}
	\min_{\epsilon_1,\ldots,\epsilon_n\in\{\pm 1\}}\norm{\sum_{i=1}^N\epsilon_iA_i}{\infty}\cdot
\end{align*}

By giving a definition of quantum discrepancy, we have introduced a notion of discrepancy for projection systems. In the given formulation, Hermitian unitary matrices play the role of (quantum) colorings. The generalization of set systems and colorings to projection systems and quantum colorings, respectively, is compatible with the general principles of quantization procedure, and it justifies using the term ``quantum'' in the denominations. Similar to the area of the combinatorial discrepancy, bounding the quantum discrepancy for general projection systems or projection systems with additional structures constitutes an important class of problems. Moreover, since the set of quantum colorings is strictly larger than the set of combinatorial colorings, the quantum discrepancy of a set system is, in general, different from its combinatorial discrepancy. 

Investigating the quantum analogues of famous problems and results in combinatorial discrepancy theory would be a potential direction for future studies. Also, computation of the quantum discrepancy for some special classes of set systems provides us with another set of interesting problems. Bounding combinatorial discrepancy would be a potential important application of quantum discrepancy. According to the probabilistic interpretation of the quantum discrepancy of set systems, this might lead to upper bounds obtained by probabilistic constructive proofs.
\subsection*{Acknowledgements}
We would like to express our great appreciation to Professor Amir Daneshgar for his valuable and constructive suggestions. Special thanks are also given to the anonymous reviewer for her/his helpful comments.

\end{document}